\newcommand\CC{{\mathbb C}}

\newcommand\cB{{\cal B}}

\newcommand\cH{{\cal H}}
\newcommand\cI{{\cal I}}

\newcommand\cl{\colon}
\newcommand\cL{{\cal L}}
\newcommand\cM{{\cal M}}

\newcommand\cO{{\cal O}}
\newcommand\cod{\mbox{cod}}

\newcommand\cQ{{\cal Q}}

\newcommand\cU{{\cal U}}
\newcommand\cV{{\cal V}}

\newcommand\e{\epsilon}

\newcommand\es{\emptyset}

\newcommand\g{\gamma}
\newcommand\G{\Gamma}

\newcommand\Gr{{\mathbb Gr}}

\newcommand\hra{\hookrightarrow}

\newcommand\la{\langle}

\newcommand\lra{\longrightarrow}

\newcommand\n{\noindent}

\newcommand\NN{{\mathbb N}}

\newcommand\om{\omega}

\newcommand\op{\oplus}
\newcommand\ot{\otimes}
\newcommand\ov{\overline}

\newcommand\PP{{\mathbb P}}

\newcommand\QQ{{\mathbb Q}}
\newcommand\ra{\rangle}

\newcommand\rk{\mbox{rk}}

\newcommand\RR{{\mathbb R}}
\newcommand\s{\sigma}

\newcommand\Si{\Sigma}

\newcommand\sm{\setminus}

\newcommand\T{\Theta}

\newcommand\tm{\times}

\newcommand\ul{\underline}

\newcommand\vf{\varphi}

\newcommand\wh{\widehat}
\newcommand\wt{\widetilde}

\newcommand\ZZ{{\mathbb Z}}

\documentclass[10pt,a4paper]{article}
\usepackage{amsthm}
\usepackage{amsmath}
\usepackage{amssymb}
\newtheorem{thm}{Theorem}[section]
\newtheorem{clm}[thm]{Claim}
\newtheorem{cnj}[thm]{Conjecture}
\newtheorem{crl}[thm]{Corollary}
\newtheorem{crt}[thm]{Condition}

\newtheorem{dfn-lmm}[thm]{Definition-Lemma}

\newtheorem{lmm}[thm]{Lemma}
\newtheorem{prp}[thm]{Proposition}

\newtheorem{rmk}[thm]{Remark}

\begin{document}
 \title{Irreducible  symplectic $4$-folds numerically
 equivalent to $(K3)^{[2]}$}
 \author{Kieran G. O'Grady\thanks{Supported by
 Cofinanziamento MURST 2002-03 and MIUR 2003-2004}\\
Universit\`a di Roma \lq\lq La Sapienza\rq\rq}
\date{July 19 2005
\vskip 4mm \small{\it Dedicato a Barbara} \vskip
1cm
\begin{abstract}
 First steps towards a classification of
 irreducible symplectic $4$-folds whose integral $2$-cohomology
 with $4$-tuple cup product is
 isomorphic to that of $(K3)^{[2]}$. We prove
 that any such $4$-fold deforms to an irreducible
 symplectic $4$-fold of
 Type A or Type B. A $4$-fold of Type A is a
 double cover
 of a (singular) sextic
 hypersurface and a $4$-fold of Type B is
  birational to a
 hypersurface of degree at most $12$. We
 conjecture that Type B $4$-folds do not exist.
 \end{abstract}}
 \maketitle
 \section{Introduction}
 \setcounter{equation}{0}
Kodaira~\cite{kod} proved that any two $K3$
surfaces are deformation equivalent. A  $K3$
surface is the same as an irreducible symplectic
$2$-fold - recall that a compact K\"ahler
manifold is irreducible symplectic if it is
simply connected and it carries a holomorphic
symplectic form spanning $H^{2,0}$
(see~\cite{beau,huy}). A classification of
higher-dimensional
 irreducible symplectic manifolds up
 to deformation equivalence appears to be out of reach
  at the moment (see~\cite{beau,huy}). We will take
 the first steps towards a solution of the
 classification problem
 for numerical $(K3)^{[2]}$'s. We
 explain our terminology: two irreducible
 symplectic manifolds $M_1,M_2$ of dimension $2n$
 are {\it numerically equivalent} if there exists
 an isomorphism of abelian groups $\psi\cl
 H^2(M_1;\ZZ)\overset{\sim}{\lra}H^2(M_2;\ZZ)$
 such that
 $\int_{M_1}\alpha^{2n}=\int_{M_2}\psi(\alpha)^{2n}$ for
 all $\alpha\in H^2(M_1;\ZZ)$. Recall~\cite{beau} that
 if $S$ is a $K3$ then $S^{[n]}$ - the
 Douady space parametrizing length-$n$
 analytic subsets of $S$ - is an irreducible symplectic
 manifold of dimension $2n$. A {\it numerical
 $(K3)^{[2]}$} is an irreducible symplectic
 $4$-fold numerically equivalent to $S^{[2]}$
 where $S$ is a $K3$.
\begin{thm}\label{mainthm1}
Let $M$ be a numerical $(K3)^{[2]}$. Then $M$ is
deformation equivalent to one of the following:
 \begin{itemize}
\item[(1)]
An irreducible symplectic $4$-fold $X$ carrying an anti-symplectic involution
 $\phi\cl X\to X$ such that the quotient
 $X/\la\phi\ra$ is isomorphic to a
 sextic hypersurface $Y\subset\PP^5$. Let
 $f\cl X\to Y$ be the quotient map and
 $H:=f^{*}\cO_{Y}(1)$; the fixed locus of $\phi$ is
a smooth irreducible Lagrangian surface
 $F$ such that
\begin{equation}\label{cherndieffe}
 c_2(F)=192,\quad \cO_F(2K_F)\cong\cO_F(6H),
   \quad c_1(F)^2=360.
\end{equation}
\item[(2)]
An irreducible symplectic $4$-fold $X$ admitting
a rational map $f\cl X\cdots>\PP^5$ which is
birational onto its image $Y$, with $6\le \deg
Y\le 12$.
\end{itemize}
\end{thm}
We give a brief outline of the proof of the theorem. By applying  surjectivity
of the period map and Huybrechts' projectivity criterion~\cite{huy,huy-err} we
will be able to deform $M$ to an irreducible symplectic $4$-fold $X$ such that
Items~(1) through~(6) of Proposition~(\ref{hodgeprop}) hold. The first item
gives (via Hirzebruch-Riemann-Roch and Kodaira Vanishing) that there is an ample
divisor $H$ on $X$ such that
\begin{equation}\label{hcond}
  \int_X c_1(H)^4=12,\quad h^0(\cO_X(H))=6.
\end{equation}
Let $h:=c_1(H)$; Items~(2), (3) and~(4) state that $h$ generates
$H^{1,1}_{\ZZ}(X)$ and that $H^4(X)$ has no rational Hodge substructures
other than those forced by $h$ and the Beauville quadratic form.
Items~(5)-(6) imply, via Proposition~(\ref{intprop}),  the following {\it
Irreducibility property of $|H|$}: if $D_1,D_2\in |H|$ are
 distinct then
$D_1\cap D_2$ is a reduced and irreducible
surface in $X$. Next we will study the rational
map $f\cl X\cdots>|H|^{\vee}\cong\PP^5$. A
straightforward argument based on ampleness of
$H$ and the Irreducibility property of $|H|$ will
show that either Item~(1) or Item~(2) of
Theorem~(\ref{mainthm1}) holds or $Y:=Im(f)$ is
one of the following
\begin{itemize}
  \item [(a)]
   a $3$-fold of degree at most $6$,
  \item [(b)]
   a $4$-fold of degree at most $4$.
\end{itemize}
We will prove that (a) or (b) cannot hold arguing
by contradiction: assuming that (a) or (b) holds
we will get that either $H^4(X)$ has a
non-existant Hodge substructure or the
Irreducibilty property of $|H|$ does not hold -
with the exception of $Y$ a normal quartic
$4$-fold, this case will require an ad hoc
argument. Thus we will need to analyze $3$-folds
and $4$-folds in $\PP^5$ of low degree. In
particular we will prove some results on cubic
$4$-folds $Y\subset\PP^5$ which might be of
independent interest. First we will show that if
$\dim(sing Y)\ge 1$ then $Y$ contains a plane.
Secondly we will prove that if $Y$ is singular
with isolated singularities and it does not
contain planes then $Gr^{W}_4 H^4(Y)$ contains a
Hodge substructure isomorphic to the
transcendental part of the $H^2$ of a $K3$
surface (shifted by $(1,1)$), namely the minimal
desingularization of the set of lines in $Y$
through any of its singular points. This result
should be equivalent to a statement about
degenerations of the variety $F(Y)$ parametrizing
lines on a cubic $4$-folds $Y\subset\PP^5$ -
recall that if $Y$ is smooth  then $F(Y)$ is a
deformation of $(K3)^{[2]}$ (see~\cite{beau-don})
and if $Y$ is singular then $F(Y)$ is
singular~\cite{joe}. The relevant statement is
the following. Let $\cU$ be the parameter space
for cubic $4$-folds $Y\subset\PP^5$ not
containing a plane: there exists  a finite cover
$\wt{\cU}\to\cU$  such that the pull-back to
$\wt{\cU}$ of the family over $\cU$ with fiber
$F(Y)$ at $[Y]$ has a simultaneous resolution of
singularities. The proof of
Theorem~(\ref{mainthm1}) should be compared to
that given in~\cite{semk3} of Kodaira's theorem
on deformation equivalence of $K3$ surfaces. The
general strategies are the same however we have
to work harder and the result is not as
conclusive as Kodaira's\footnote{We do have a
proof that $f$ cannot be birational onto its
image with $\deg(Im(f))\le 8$.}.
\begin{cnj}\label{doppio}
Suppose that $X$ is a numerical $(K3)^{[2]}$ and
that Items~(1) through~(6) of
Proposition~(\ref{hodgeprop}) hold. Then Item~(1)
of Theorem~(\ref{mainthm1}) holds.
\end{cnj}
We notice that if $X$ satisfies Item~(1) of
Theorem~(\ref{mainthm1}) then any small
deformation of $X$ that keeps $c_1(H)$ of type
$(1,1)$ is a variety which again satisfies
Item~(1) and the hyperplane class on the deformed
variety is the deformation of the hyperplane
class on $X$ - see Proposition~(\ref{stabile}).
In other words the conjecture is stable for small
deformations. If the above conjecture is true
then  any numerical $(K3)^{[2]}$ is deformation
equivalent to an $X$ as in Item~(1) of
Theorem~(\ref{mainthm1}). In another paper we
prove that the quotient $Y$ belongs to the set of
sextic hypersurfaces described by
Eisenbud-Popescu-Walter in Example~(9.3)
of~\cite{epw}. We will also show that if $Z$ is a
generic EPW-sextic and $W\to Z$ is the natural
double cover then $W$ is deformation of
$(K3)^{[2]}$; this will imply that if
Conjecture~(\ref{doppio}) holds then any
numerical $(K3)^{[2]}$ is a deformation of
$(K3)^{[2]}$.

 \vskip 3mm
\n
 {\bf Notation:} If $X$ is a topological space then
$H^{*}(X)$ denotes cohomology with complex
coefficients.

Topology of algebraic varieties (or analytic
spaces) will be either the classical topology or
the Zariski topology: in general it will be clear
from the context in which topology we are
working.

Let $X$ be a smooth projective variety. If $W$ is
a closed subscheme of $X$ of pure dimension $d$
we let
\begin{equation}\label{fundcyc}
  [W]\in Z_d(X)
\end{equation}
be the fundamental cycle associated to $W$ as
in~\cite{ful}, p.~15.

Let $\PP(V)$ be a projective space. If
$A\subset\PP(V)$ we let $span(A)\subset\PP(V)$ be
the span of $A$, i.e.~the intersection of all
linear subspaces containing $A$. If
$A,B\subset\PP(V)$ we let
\begin{equation}\label{joindef}
  J(A,B):=\bigcup\limits_{p\in A,q\in B}span(p,q).
\end{equation}
If $A,B$ are closed and $A\cap B=\es$ then $J(A,B)$ is closed - in general
$J(A,B)$ is not closed.

Let $X$ be a scheme and $x\in X$ a (closed)
point; we let $\T_x X$ be the Zariski tangent
space to $X$ at $x$. Now assume that $X$ is a
subscheme of a projective space $\PP(V)$. Then
$\T_x X\subset\T_x\PP(V)$: the {\it projective
tangent space to $X$ at $x$} is the unique linear
subspace
\begin{equation}\label{projtan}
  T_x X\subset\PP(V)
\end{equation}
containing $x$ whose Zariski tangent space at $x$ is equal to $\T_x X$.
\vskip 3mm

 \n
 {\bf Acknowledgements:} Initially I proved the
results of Sections~(\ref{depardieu})-(\ref{dimomainthm}) for $X$ a deformation
of $(K3)^{[2]}$ provided with an ample divisor $H$ of square $2$ for the
Beauville quadratic form. Claire Voisin observed that the proofs had to be valid
for symplectic $4$-folds satisfying suitable cohomological hypotheses: I thank
Claire for her precious observation.
\section{Preliminaries}
\subsection{Beauville's form and Fujiki's constant}
\label{beaufuj}
 \setcounter{equation}{0}
Let $M$ be an irreducible symplectic manifold of
dimension $2n$. By Beauville and Fujiki
(see~\cite{beau} and Thm.~(4.7)
 of~\cite{fuj}) there exist a rational positive
number $c_M$ and an integral indivisible
 symmetric bilinear form $(,)_M$ on $H^2(M)$
characterized by the following properties. First
 $(,)_M$ is positive definite on the span of
 $\{\s+\ov{\s}\}_{\s\in H^{2,0}(M)}$ and an
 arbitrarily
  chosen K\"ahler
class. Secondly we have the equality
\begin{equation}\label{fujrel}
 \int_M \alpha^{2n}=c_M(\alpha,\alpha)_M^n,\qquad \alpha
 \in H^2(M).
\end{equation}
 Thus $H^2(M;\ZZ)$ has a canonical structure of
lattice.  If two irreducible symplectic manifolds
of the same dimension have the same Beauville
form and Fujiki constant then by~(\ref{fujrel})
they are numerically equivalent. The converse is
\lq\lq almost true\rq\rq. In fact let $\om\in
H^{1,1}(M;\RR)$ be a K\"ahler class;
by~(\ref{fujrel}) the primitive (with respect to
$\om$) cohomology $H^2(M)_{prim}$ is equal to
$\om^{\bot}$ (orthogonality is with respect to
$(,)_M$) and hence by
 the Hodge index Theorem the signature of $(,)_M$ is
  $(3,b_2(M)-3)$. It follows that
  if two irreducible symplectic
manifolds $M_1,M_2$ of dimension $2n$ are
numerically equivalent then they have the same
Beauville form and Fujiki constant unless $n$ is
even and $b_2(M_1)=b_2(M_2)=6$: in this case
numerical equivalence implies that $(,)_{M_1}=\pm
(,)_{M_2}$ (and $c_{M_1}=c_{M_2}$).
 Let $\Lambda$ be the lattice given by
\begin{equation}\label{lambdadef}
\Lambda:=U^{\op 3}\op (-E_8)^{\op 2}\op (-2),
\end{equation}
where $U$ is the standard hyperbolic plane. Let
$S$ be a $K3$ surface; the Beauville form and
Fujiki constant of $S^{[2]}$ are given
(see~\cite{beau}) by

\begin{equation}\label{hilbcase}
H^2(S^{[2]};\ZZ)\cong \Lambda,\qquad
c_{S^{[2]}}=3.
\end{equation}
Thus a numerical $(K3)^{[2]}$ is an irreducible
symplectic $4$-fold $M$ such that
\begin{equation}\label{numhyp}
 H^2(M;\ZZ)\cong \Lambda,\qquad
c_M=3.
\end{equation}
In particular $b_2(M)=23$; as is well-known -
see~\cite{guan,verbit} - this implies that
\begin{equation}\label{symmeq}
H^3(M;\QQ)=0,\qquad Sym^2 H^2(M;\QQ)\overset{\sim}{\lra} H^4(M;\QQ),
\end{equation}
where the second isomorphism is given by
cup-product. The equations of~(\ref{symmeq}) will
be crucial for what follows.
\subsection{Quadratic forms on $V$ and $S^2V$}
\label{determinante}
 \setcounter{equation}{0}
Let $A$ be a ring and $V$ be an $A$-module. Let
$(V\ot V)^{+},(V\ot V)^{-}\subset V\ot V$ be the
submodules of tensors which are invariant,
 respectively
anti-invariant, for the involution of $V\ot V$ interchanging the factors. We let
$Sym_2 V:=(V\ot V)^{+}$ and $Sym^2 V:=V\ot V/(V\ot V)^{-}$.
 Assume that $(,)$ is
a symmetric bilinear form on $V$; we let $\la,\ra$ be the unique symmetric
bilinear form on $S^2 V$ such that
\begin{equation}\label{donpol}
    \langle\alpha_1\alpha_2,\alpha_3\alpha_4\rangle=
    (\alpha_1,\alpha_2)(\alpha_3,\alpha_4)+
    (\alpha_1,\alpha_3)(\alpha_2,\alpha_4)
    +(\alpha_1,\alpha_4)(\alpha_2,\alpha_3)
\end{equation}
for $\alpha_1,\ldots,\alpha_4\in V$.
Using~(\ref{fujrel}) and the second equality
of~(\ref{numhyp}) we get the following.
\begin{rmk}\label{collegamento}
Let $M$ be a numerical $(K3)^{[2]}$. The
intersection form on
\begin{equation}\label{juve}
  Sym^2 H^2(M)\cong H^4(M)
\end{equation}
is the bilinear form constructed as above from
$V:=H^2(M)$ and $(,):=(,)_M$.
\end{rmk}
\section{The deformation}
\setcounter{equation}{0}
Let $M$ be a numerical $(K3)^{[2]}$. We will show
that $M$ can be deformed to a projective
irreducible symplectic $4$-fold $X$ such that
$H^{*}(X)$ has few integral Hodge substructure.
First we  introduce the tautological rational
Hodge substructures of $H^{*}(X)$ for $X$ a
numerical $(K3)^{[2]}$ with an $h\in
H^{1,1}_{\QQ}(X)$ such that
\begin{equation}\label{integrnon}
  \int_X h^4\not=0.
\end{equation}
To simplify notation we let $(,)$ be the
Beauville form of $X$; thus~(\ref{integrnon}) is
equivalent by~(\ref{fujrel}) to $(h,h)\not=0$. We
have an orthogonal direct sum decomposition
\begin{equation}\label{h2decomp}
 H^2(X;\CC)=\CC h\op_{\bot} h^{\bot}
\end{equation}
into Hodge substructures of levels $0$ and $2$ respectively.
By~(\ref{symmeq}) we have a direct sum decomposition
\begin{equation}\label{symmdecomp}
    H^4(X;\CC)=\CC h^2\op \left(\CC h\ot h^{\bot}\right)\op Sym^2(h^{\bot})
\end{equation}
into Hodge substructures of levels $0$, $2$ and $4$
 respectively.   There is a refinement of
Decomposition~(\ref{symmdecomp}); to
 explain this we need to introduce the dual
of Beauville's form. Let
$q\in Sym_2\left(H^2(X)^{\vee}\right)$ be Beauville's
symmetric bilinear form; it is non-degenerate~\cite{beau}
 and hence it defines
an isomorphism $L_q\cl
H^2(X)\overset{\sim}{\to}H^2(X)^{\vee}$. Let
$\Pi_2\cl Sym_2 H^2(X)\to Sym^2 H^2(X)$ be the
composition of the inclusion $Sym_2 H^2(X)\hra
H^2(X)\ot H^2(X)$ and the projection map
$H^2(X)\ot H^2(X)\to Sym^2 H^2(X)$. Let
\begin{equation}\label{qdual}
    q^{\vee}:=\Pi_2\circ Sym_2(L_q^{-1})(q)\in Sym^2 H^2(X).
\end{equation}
Explicitly: let $\{\alpha_1,\ldots,\alpha_{23}\}$
be a basis of $H^2(X)$ and
$\{\alpha^{\vee}_1,\ldots,\alpha^{\vee}_{23}\}$
be the dual basis. Thus
\begin{equation}\label{explq}
    q=\sum_{ij}g_{ij}\alpha^{\vee}_i\ot \alpha^{\vee}_j
\end{equation}
where $(g_{ij})$ is a symmetric matrix. Then
\begin{equation}\label{explqdual}
    q^{\vee}=\sum_{ij}m_{ij}\alpha_i\alpha_j,
    \qquad (m_{ij})=(g_{ij})^{-1}.
\end{equation}
We know that $q$ is integral and that
$(\alpha_1,\alpha_2)=0$ if $\alpha_i\in
H^{r_i,2-r_i}(X)$ with $r_1+r_2\not=2$; this
implies that
\begin{equation}\label{qpure}
 q^{\vee}\in H^{2,2}_{\QQ}(X).
\end{equation}
In terms of Decomposition~(\ref{symmdecomp}) we
have $q^{\vee}\in\CC h^2\op Sym^2(h^{\bot})$.
More precisely let $q_h:=q|_{h^{\bot}}$ and let
$q_h^{\vee}\in Sym^2(h^{\bot})$ be its dual (this
makes sense because $(h,h)\not=0$ and hence $q_h$
is non-degenerate); then
\begin{equation}\label{qdecomp}
 q^{\vee}=(h,h)^{-1}h^2+q_h^{\vee}.
\end{equation}
Let $\langle,\rangle$ be the intersection form on
$H^4(X)$ - the notation is consistent with  that
of Subsection~(\ref{determinante}), see
Remark~(\ref{collegamento}) - and let
\begin{equation}\label{wdef}
    W(h):=(q^{\vee})^{\bot}\cap Sym^2(h^{\bot}),
\end{equation}
where the first orthogonality is with respect to
$\la,\ra$ and the second is with respect to
$(,)$.
\begin{clm}\label{wnontriv}
Keeping notation as above, $W(h)$ is a codimension-$1$ rational
sub Hodge structure of $Sym^2(h^{\bot})$, and we have a direct sum
decomposition
\begin{equation}\label{smalldecomp}
\CC h^2\op Sym^2(h^{\bot})=\CC h^2\op \CC q^{\vee}\op W(h).
\end{equation}
\end{clm}
\begin{proof}
$W(h)$ is a sub Hodge structure because $q^{\vee}$ is rational of type $(2,2)$;
let's show that
\begin{equation}\label{notinc}
 Sym^2(h^{\bot})\not\subset (q^{\vee})^{\bot}.
\end{equation}
From Remark~(\ref{collegamento})  one gets that
\begin{equation}\label{qdualint}
 \langle q^{\vee},\alpha\beta\rangle=
 25(\alpha,\beta),\qquad \alpha,\beta\in H^2(X).
\end{equation}
From this we get immediately~(\ref{notinc}) and
thus $W(h)$ has codimension $1$. Now let's prove
that we have~(\ref{smalldecomp}).
By~(\ref{qdecomp}) $h^2$ and $q^{\vee}$ are
linearly independent and hence it suffices to
show that
\begin{equation}\label{heydude}
\left(\CC h^2\op \CC q^{\vee}\right)\cap W(h)=\{0\}.
\end{equation}
It follows from~(\ref{qdualint}) that
\begin{equation}\label{qdualqdual}
 \langle q^{\vee}, q^{\vee}\rangle=25\cdot 23,
\end{equation}
and hence
\begin{equation}\label{cucs}
  \left(\CC h^2\op \CC
q^{\vee}\right)\cap(q^{\vee})^{\bot}=
\CC(23h^2-(h,h)q^{\vee}).
\end{equation}
On the other hand by~(\ref{qdecomp}) we have
\begin{equation}\label{spqr}
  \left(\CC h^2\op \CC q^{\vee}\right)\cap
Sym^2(h^{\bot})=\CC(h^2-(h,h)q^{\vee}).
\end{equation}
Equation~(\ref{heydude}) follows immediately
from~(\ref{cucs})-(\ref{spqr}).
\end{proof}
By the above claim  we have a  decomposition
\begin{equation}\label{newsymmdecomp}
    H^4(X;\CC)=\left(\CC h^2\op \CC q^{\vee}\right)\op\left(\CC h\ot h^{\bot}\right)\op W(h)
\end{equation}
into sub-H.S.'s of levels $0$, $2$ and $4$
respectively. The following is the main result of
this section.
\begin{prp}\label{hodgeprop}
Keep notation as above. Let $M$ be a numerical
$(K3)^{[2]}$. There exists an irreducible
symplectic manifold $X$ deformation equivalent to
 $M$ such that:
 \begin{itemize}
\item[(1)]
 $X$ has an ample divisor $H$  with  $(h,h)=2$,
 where $h:=c_1(H)$,
\item[(2)]
 $H^{1,1}_{\ZZ}(X)=\ZZ h$,
\item[(3)]
Let $\Si\in Z_1(X)$ be an integral algebraic
$1$-cycle on $X$ and $cl(\Si)\in
H^{3,3}_{\QQ}(X)$ be its Poincar\'e dual. Then
$cl(\Si)=m h^3/6$ for some $m\in\ZZ$.
\item[(4)]
if $V\subset H^4(X)$ is a rational sub Hodge
structure then $V=V_1\op V_2\op V_3$ where
$V_1\subset\left(\CC h^2\op \CC q^{\vee}\right)$,
$V_2$ is either $0$ or equal to $\CC h\ot
h^{\bot}$ and $V_3$ is either $0$ or equal to
$W$.
\item[(5)]
the image of $h^2$ in $H^4(X;\ZZ)/Tors$ is indivisible,
\item[(6)]
$H^{2,2}_{\ZZ}(X)/Tors\subset\ZZ(h^2/2) \op\ZZ
(q^{\vee}/5)$.
\end{itemize}
\end{prp}
The proof of the proposition will be given after
some preliminary results. We recall Huybrechts'
Theorem on surjectivity of the global period
map~\cite{huy,huy-err} - in the context of
numerical $(K3)^{[2]}$'s. Let $M$ be a numerical
$(K3)^{[2]}$ and $\cM$ be the moduli space of
marked irreducible symplectic manifolds
deformation equivalent to $M$; thus a point of
$\cM$ is an equivalence class of couples
$(X,\psi)$ where $X$ is an irreducible symplectic
manifold deformation equivalent to $M$ and
$\psi\cl \Lambda \overset{\sim}{\lra} H^2(X;\ZZ)$
is an isometry of lattices ($\Lambda$ is the
lattice~(\ref{lambdadef})). The couples
$(X,\psi)$ and $(X',\psi')$ are equivalent if
there exists an isomorphism $f\cl X\to X'$ such
that $H^2(f)\circ\psi' =\pm\psi$. If $t\in\cM$ we
let $(X_t,\psi_t)$ be a representative of $t$. It
is known that $\cM$ is a non-separated complex
analytic space, see Thm.(2.4) of~\cite{lp}. The
period domain $Q\subset\PP(\Lambda\ot \CC)$ is
given by
\begin{equation}\label{perdom}
    Q:=\{[\s]\in \PP(\Lambda\ot \CC)|\ (\s,\s)_{\Lambda}=0.
    \ (\s,\ov{\s})_{\Lambda}>0\}
\end{equation}
where $(,)_{\Lambda}$ is the symmetric bilinear
form on $\Lambda$. The period map is given by
$$\begin{array}{ccc}
  \cM & \overset{P}{\lra} & Q \\
  (X,\psi) & \mapsto & \psi^{-1}(H^{2,0}(X)). \\
\end{array}$$
Here and in the following $\psi$ denotes both the
isometry $\Lambda \overset{\sim}{\lra}
H^2(X;\ZZ)$ and its linear extension
$\Lambda\ot\CC\to H^2(X;\CC)$. The map $P$ is
locally an isomorphism, see~\cite{beau}). Let
$\cM^0$ be a connected component of $\cM$.
Huybrechts' Theorem on surjectivity of the global
period map (Thm.~(8.1) of~\cite{huy}) states that
the restriction of $P$ to $\cM^0$ is surjective.
Given $\alpha\in \Lambda$ we let
\begin{equation}\label{emmealfa}
  \cM^0_{\alpha}:=\{t\in \cM^0|\ \psi_t(\alpha)
    \in H^{1,1}_{\ZZ}(X_t)\}.
\end{equation}
\begin{lmm}\label{primhs}
Let $M$ be a numerical $(K3)^{[2]}$ and $\cM$ be
the moduli space of marked irreducible symplectic
manifolds deformation equivalent to $M$. Let
$\alpha\in \Lambda$ with $(\alpha,\alpha)\not=0$.
For $t\in\cM^0_{\alpha}$ outside of a countable
union of proper analytic subsets we have:
 \begin{itemize}
\item[(1)]
 $H^{1,1}_{\QQ}(X_t)=\QQ \psi_t(\alpha)$,
\item[(2)]
any rational sub Hodge structure of
$W(\psi_t(\alpha))$  is trivial.
\end{itemize}
\end{lmm}
\begin{proof}
Let
\begin{equation}\label{lgdef}
    L_{\alpha}:=\{[\s]\in Q|\ (\s,\alpha)_{\Lambda}=0\}.
\end{equation}
As is easily checked $L_{\alpha}$ is a non-empty
codimension $1$ subvariety of $Q$ and furthermore
\begin{equation}\label{invlg}
   \cM^0_{\alpha}=P^{-1}(L_{\alpha}).
\end{equation}
By surjectivity of the period map
$\cM^0_{\alpha}$ is non-empty of dimension $20$.
 It is well-known that the set of $t\in
\cM^0_{\alpha}$ for which~(1) does not hold is a
countable union of proper analytic subsets of
$\cM^0_{\alpha}$(see~\cite{huy}).  Next we show
that the set of $t\in\cM^0_{\alpha}$ for which
(2) does not hold is also a
 countable union of proper
analytic subsets of $\cM^0_{\alpha}$. Let
$$W(\alpha):=Sym^2(\alpha^{\bot})\cap
(q_{\Lambda}^{\vee})^{\bot}\subset
Sym^2(\Lambda\ot\CC)$$
where $q_{\Lambda}$ is the quadratic form on
$\Lambda$. For a linear subspace $V\subset
W(\alpha)$ defined over $\QQ$ let
\begin{equation}\label{cmvdef}
    \cM^0_{\alpha}(V):=\{t\in\cM^0_{\alpha}|\
    \mbox{$Sym^2(\psi_t)(V)$ is a sub-H.S.~of
    $W(\psi_t(\alpha))$}\}.
\end{equation}
Since the set of subspaces $V\subset W(\alpha)$
defined over $\QQ$ is countable it suffices to
prove that $\cM^0_{\alpha}(V)$ is a proper
analytic subset of $\cM^0_{\alpha}$ whenever
$V\not=0$ or $V\not= W(\alpha)$. It is well-known
that $\cM^0_{\alpha}(V)$ is an analytic subset of
$\cM^0_{\alpha}$. Assume that $\cM^0_{\alpha}(V)$
contains a non-empty open subset $U\subset
\cM^0_{\alpha}$: we will show that either $V=0$
or $V=W(\alpha)$. We have
\begin{itemize}
\item[(a)]
 $Sym^2(\psi_t)(V)\cap H^{4,0}(X_t)\not=\{0\}$
 for all $t\in U$, or
\item[(b)] $Sym^2(\psi_t)(V)\cap H^{4,0}(X_t)=\{0\}$
for all $t\in U$.
\end{itemize}
Assume that (a) holds. Then
$$Sym^2(\psi_t)(V)\supset H^{4,0}(X_t)=
H^{2,0}(X_t)\wedge H^{2,0}(X_t)$$
for all $t\in U$ and hence
\begin{equation}\label{pinodaniele}
  V\supset\{\s^2|\ [\s]\in P(U)\},
\end{equation}
where $P$ is the period map. Let $\cV\subset
L_{\alpha}$ be open and non-empty: as is easily
verified
\begin{equation}\label{woody}
   span\{[\s^2]|\ [\s]\in \cV\}
=\PP(W(\alpha)).
\end{equation}
Since $P(U)$ is an open non-empty subset of
$L_{\alpha}$ we get by~(\ref{pinodaniele}) that
$V=W(\alpha)$. Now assume that (b) holds. Then
$$\langle Sym^2(\psi_t)(V),H^{0,4}(X_t)\rangle=0$$
for all $t\in U$ and hence $V\bot\{\ov{\s}^2|\
[\s]\in P(U)\}$. Arguing as above we get that
$V=\{0\}$.
\end{proof}
We will apply Lemma~(\ref{primhs}) with a
particular choice of $\alpha$. First we prove two
preliminary results.
\begin{lmm}\label{tuttiequiv}
The vectors $\alpha\in\Lambda$ with
\begin{equation}\label{squaretwo}
  (\alpha,\alpha)_{\Lambda}=2
\end{equation}
 belong to a single $O(\Lambda)$-orbit.
\end{lmm}
\begin{proof}
Let
$$\wt{\Lambda}:=U^{\op 3}\op (-E_8)^{\op 2}\op U.$$
Choose an embedding $\Lambda\subset\wt{\Lambda}$
such that $\Lambda^{\bot}=\ZZ\g$ where $\g\in U$
is a vector with $(\g,\g)_{\Lambda}=2$. Given
$\alpha_1,\alpha_2\in\Lambda$ with
$(\alpha_i,\alpha_i)_{\Lambda}=2$ the lattices
$\ZZ\g\op\ZZ\alpha_1$ and $\ZZ\g\op\ZZ\alpha_2$
are saturated and isometric. By a standard result
on lattices (see Theorem~1, p.~578
of~\cite{shafa}) there exists $g\in
O(\wt{\Lambda})$ with $g(\g)=\g$ and
$g(\alpha_1)=\alpha_2$. Since $g$ sends
$\Lambda=\g^{\bot}$ to itself it restricts to an
isometry of $\Lambda$ taking $\alpha_1$ to
$\alpha_2$.
\end{proof}
\begin{lmm}\label{alphachoice}
Let $M$ be a numerical $(K3)^{[2]}$. Let $\cM$ be
the moduli space of marked irreducible symplectic
manifolds deformation equivalent to $M$ and let
$\cM^0$ be a connected component of $\cM$.
Suppose that $\alpha_1,\alpha_2\in\Lambda$
satisfy
\begin{equation}\label{intermatrix}
(\alpha_1,\alpha_1)_{\Lambda}=(\alpha_2,\alpha_2)_{\Lambda}=2,\quad
(\alpha_1,\alpha_2)_{\Lambda}\equiv 1\mod{2}.
\end{equation}
There exists $1\le i\le 2$  such that for every
 $t\in\cM^0$ the class of $\psi_t(\alpha_i)^2$ in
 $H^4(X_t;\ZZ)/Tors$ is indivisible.
\end{lmm}
\begin{proof}
First notice that it suffices to show that for
one $t_0\in\cM^0$ there exists $1\le i\le 2$ such
that $\psi_{t_0}(\alpha_i)^2$  is indivisible; in
fact for any other $t\in\cM^0$ there exists a
diffeomorphism $f\cl
X_{t_0}\overset{\sim}{\lra}X_{t}$ such that
$H^2(f)\circ \psi_{t}=\psi_{t_0}$ and hence if
$\psi_{t}(\alpha_i)^2$ is divisible then
$\psi_{t_0}(\alpha_i)^2$ is divisible too. Next
we claim that for $i=1,2$ the class of
$\psi_t(\alpha_i)^2$  is divisible at most by
$2$. First notice that there exists
$\beta_i\in\Lambda$ with
\begin{equation}\label{betachoice}
(\alpha_i,\beta_i)_{\Lambda}=1,\quad
(\beta_i,\beta_i)_{\Lambda}=0.
\end{equation}
In fact  by Lemma~(\ref{tuttiequiv}) it suffices
to exhibit $\alpha',\beta'\in\Lambda$  such that
\begin{equation}\label{eqprimo}
  (\alpha',\alpha')_{\Lambda}=2,\quad (\alpha',\beta')_{\Lambda}=1,\quad
(\beta',\beta')_{\Lambda}=0,
\end{equation}
and this is a trivial exercise. Now let $\beta_i$
be as above. Then
\begin{equation}\label{piccolo}
  \langle \psi_t(\alpha_i)^2,\psi_t(\beta_i)^2 \rangle=2
\end{equation}
and this proves that $\psi_t(\alpha_i)^2$ is
divisible at most by $2$. Now we prove the lemma
arguing by contradiction. Assume that
$\psi_t(\alpha_i)^2$ is divisible by $2$ (modulo
torsion) for $i=1$ and $i=2$; thus
\begin{equation}\label{arbore}
  \psi_t(\alpha_i)^2=2\g_i+\xi_i
\end{equation}
where $\g_i\in H^4(M;\ZZ)$ and $\xi\in
Tors(H^4(M;\ZZ))$. By Remark~(\ref{collegamento})
we have
\begin{equation}\label{bellaciao}
  \langle \psi_t(\alpha_1)^2, \psi_t(\alpha_2)^2
  \rangle=6.
\end{equation}
On the other hand by~(\ref{arbore}) the left-hand
side is equal to $4\la \g_1,\g_2\ra$,
contradiction.
\end{proof}
\n
 {\bf Proof of Proposition~(\ref{hodgeprop}).}
Let $\cM$ be the moduli space of marked
irreducible symplectic manifolds deformation
equivalent to $M$ and let $\cM^0$ be a connected
component of $\cM$. By Lemma~(\ref{alphachoice})
there exists $\alpha\in\Lambda$ with
$(\alpha,\alpha)=2$ such that for every
$t\in\cM^0$ the class of $\psi_t(\alpha)^2$ in
$H^4(X_t;\ZZ)/Tors$ is indivisible. Let
$t\in\cM^0_{\alpha}$ satisfying Items~(1)-(2) of
Lemma~(\ref{primhs}). Set $X:=X_t$. Since
$\psi_t(\alpha)\in H^{1,1}_{\ZZ}(X)$ and
$(\psi_t(\alpha),\psi_t(\alpha))=2$ we know that
$X$ is projective by Huybrechts' projectivity
criterion~\cite{huy}: since
$H^{1,1}_{\ZZ}(X)=\ZZ\psi_t(\alpha)$ either
$\psi_t(\alpha)$ or $-\psi_t(\alpha)$ is the
class of an ample divisor. Let
$h:=\psi_t(\alpha)$ in the former case and
$h:=-\psi_t(\alpha)$ in the latter case. We let
$H$ be a divisor with $c_1(H)=h$. Let's prove
that Items~(1)-(5) of
Proposition~(\ref{hodgeprop}) hold for $(X,H)$.
Of course $X$ is a deformation of $M$ by
definition. (1)-(2): They hold by construction.
(3): By Item~(2) and Hard Lefschetz we have
$H^{3,3}_{\QQ}(X)=\QQ h^3$ and hence $cl(\G)=x
h^3$ for some $x\in\QQ$. There exists $e\in
H^2(X;\ZZ)$ with $(e,h)=1$,
see~(\ref{betachoice}), and hence
$$\ZZ\ni \int_{\G}e=\la
xh^3,e\ra=3x(h,h)(h,e)=6x.$$
(4): Follows from Item~(2) of
Lemma~(\ref{primhs}), from the fact that $\CC
h\ot h^{\bot}$ has no non-trivial sub-H.S.'s and
an easy argument based on the observation that
the three summands of
Decomposition~(\ref{newsymmdecomp}) have pairwise
distinct levels. (5): Holds by our choice of
$\alpha$, thanks to Lemma~(\ref{alphachoice}).
(6): First we show that
\begin{equation}\label{c2formula}
    c_2(X)=6q^{\vee}/5 \quad \mbox{in $H^4(X;\QQ)$.}
\end{equation}
It is well-known that any $\theta\in Sym^2 H^2(X;\QQ)\cap H^{2,2}(X)$ which
stays of type $(2,2)$ for all deformations of $X$ is a multiple of $q^{\vee}$:
to prove it let $u\in\cM$ and use $L_{q_u}\cl
H^2(X_u)\overset{\sim}{\lra}H^2(X_u)^{\vee}$ to produce from $\theta$ a
$\theta^{'}_u\in Sym^2 H^2(X;\QQ)^{\vee}$ with $\theta^{'}_u(\s_u,\s_u)=0$ for
$\s_u\in H^{2,0}(X_u)$. Applying this to $\theta=c_2(X)$ we get that $c_2(X)=a
q^{\vee}$ for some $a\in\QQ$. Applying  Hirzebruch-Riemann-Roch and keeping in
mind that all odd Chern classes of $X$ vanish we get that
\begin{equation}\label{toddofx}
    3=\chi(\cO_X)=\frac{1}{240}\left(c_2(X)^2-\frac{1}{3}c_4(X)\right).
\end{equation}
By~(\ref{symmeq}) we know that
\begin{equation}\label{ciquattro}
  c_4(X)=324
\end{equation}
and hence it follows that
\begin{equation}\label{c2quadro}
  c_2(X)^2=828.
\end{equation}
Applying Formula~(\ref{qdualqdual}) we get that
$a=\pm 6/5$. On the other hand Theorem~(1.1)
of~\cite{miyaoka} together with~(\ref{qdualint})
gives that
\begin{equation}\label{}
  0\le \la c_2(x),h^2\ra =
  \la a q^{\vee},h^2\ra = 50 a.
\end{equation}
This proves~(\ref{c2formula}).  Since
$2q^{\vee}\in Sym^2 H^2(X;\ZZ)$
Formula~(\ref{c2formula}) gives that
\begin{equation}\label{c2over3}
    H^4(X;\ZZ)/Tors\ni (2c_2(X)-2q^{\vee})=2q^{\vee}/5
    =c_2(X)/3.
\end{equation}
In particular
\begin{equation}\label{22includes}
    \Omega(h):=\ZZ h^2\op\ZZ(2q^{\vee}/5)\subset (H^{2,2}_{\ZZ}(X)/Tors)
\end{equation}
By Item~(4) of the proposition we know that
$h^{2,2}_{\QQ}=2$ and hence $\Omega(h)$ is of
finite index in $H^{2,2}_{\ZZ}(X)/Tors$. A
straightforward computation
(use~(\ref{qdualqdual}) and~(\ref{qdualint}))
shows that
\begin{equation}\label{smalldisc}
    \mbox{discr}\left(\langle,\rangle|_{\Omega(h)} \right)=2^6\cdot 11,
\end{equation}
and hence
\begin{equation}\label{smallindex}
[H^{2,2}_{\ZZ}(X)/Tors:\Omega(h)]\le 8.
\end{equation}
Now let $xh^2+y(2q^{\vee}/5)\in
H^{2,2}_{\ZZ}(X)/Tors$: we must show that
$2x\in\ZZ$ and $2y\in\ZZ$. Let $\beta\in
H^2(X;\ZZ)$ with $(h,\beta)=1$ and
$(\beta,\beta)=0$: such a $\beta$ exists, see the
proof of Lemma~(\ref{alphachoice}).
Using~(\ref{qdualint}) we get that
$$\ZZ\ni\langle xh^2+y(2q^{\vee}/5),\beta^2\rangle=2x.$$
Next let $\g,\delta\in H^2(X;\ZZ)$ with
$(\g,\delta)=1$. Then
$$\ZZ\ni\langle xh^2+y(2q^{\vee}/5),\g\delta\rangle=2x(1+(h,\g)
(h,\delta))+10y.$$
Since $2x\in\ZZ$ we get that $10y\in\ZZ$. By~(\ref{smallindex}) we know that
$8y\in\ZZ$ and hence $2y\in\ZZ$. This finishes the proof of
Proposition~(\ref{hodgeprop}).
\begin{rmk}\label{localenough}
In the proof of Proposition~(\ref{hodgeprop}) we appealed to Huybrechts' Global
Surjectivity Theorem. It is plausible that local surjectivity is sufficient.
\end{rmk}
 If $M$ is a numerical $(K3)^{[2]}$ we cannot
 exclude the existence of a $\gamma\in H^2(M;\ZZ)$
 such that $(\gamma,\gamma)=2$ and the image of
 $\gamma^2$ in $H^4(M;\ZZ)/Tors$ is divisible by
 $2$; we only proved that it is impossible
 that $\gamma^2$  is divisible by
 $2$ for all $\gamma$ with $(\gamma,\gamma)=2$.
 If $M$ is a deformation of $(K3)^{[2]}$ the
 picture is simpler.
 \begin{prp}\label{}
 Let $M$ be a deformation of $(K3)^{[2]}$ and
  $\gamma\in H^2(M;\ZZ)$
 such that $(\gamma,\gamma)=2$. The image of
 $\gamma^2$ in $H^4(M;\ZZ)/Tors$ is not
 divisible.
 \end{prp}
 \begin{proof}
 Let $S$ be a $K3$ surface. We may assume that
 $\gamma\in H^2(S^{[2]};\ZZ)$. Let $\Delta\subset
  S^{[2]}$
 be the codimension-$1$ locus parametrizing
 non-reduced subschemes of $S$. There exists
 $\xi\in H^2(S^{[2]};\ZZ)$ such that
 $2\xi=c_1(\Delta)$. There is an orthogonal
 direct sum decomposition (see Prop.~6, p.~768
   and pp.~777-778 of~\cite{beau})
\begin{equation}\label{accaduedecomp}
 H^2(S^{[2]};\ZZ)=\mu(H^2(S;\ZZ))\oplus_{\bot}\ZZ\xi
\end{equation}
where $\mu\colon H^2(S;\ZZ)\to H^2(S^{[2]};\ZZ)$
 is the symmetrization map (Donaldson map).
 If $C\subset S$ is an algebraic curve a
 representative of $\mu(C)$ is the divisor
\begin{equation}\label{}
  \Sigma_C:=\{[Z]\in S^{[2]}|\ Z\cap
  C\not=\emptyset\}.
\end{equation}
By~(\ref{accaduedecomp}) we have
$\gamma=\mu(\alpha')-x\xi$. We know that $\gamma$
is at most divisible by $2$, see the proof of
Lemma~(\ref{alphachoice}), and hence we may add
to $\gamma$ arbitrary elements of
$2H^2(S^{[2]};\ZZ)$. Thus we may assume that
$\gamma=\mu(\pm\alpha)-\xi$ where
$(\alpha,\alpha)=4$. We can deform the complex
structure of $S$ so that either $\alpha$ or
$-\alpha$ is represented by a very ample divisor
on $S$ giving an embedding $S\subset\PP^3$. We
can furthermore assume that $S$ contains a conic
$C$. Now consider the map
\begin{equation}\label{}
   \begin{matrix}
 S^{[2]} & \overset{g}{\lra} & {\bf Gr}(1,\PP^3) \\
 [Z] & \mapsto & \langle Z\rangle
   \end{matrix}
\end{equation}
where $\langle Z\rangle $ is the line spanned by
$Z$. Let $p\colon {\bf Gr}(1,\PP^3)
 \overset{p}{\hra}\PP^5$ be the Pl\"ucker
 emebedding. Then $c_1((pg)^{*}\cO_{\PP^5}(1))=\gamma$,
 see Formula~(4.1.9) of~\cite{oginv}. Now consider the
surface $C^{(2)}\subset S^{[2]}$. Since $g$ maps
$C^{(2)}$ isomorphically onto a linear $\PP^2$ in
$\PP^5$ we get that
\begin{equation}\label{}
  \int_{C^{(2)}}\gamma^2=1.
\end{equation}
Thus the image of $\gamma^2$ in
$H^4(S^{[2]};\ZZ)/Tors$ is not divisible.
 \end{proof}
\section{The linear system $|H|$}\label{depardieu}
\setcounter{equation}{0}
Let $X,H$ be as in Proposition~(\ref{hodgeprop}).
In this section we will prove
 some basic properties
of the complete linear system $|H|$.  A key
result  is the following.
\begin{prp}\label{intprop}
Keep notation as above.
\begin{itemize}
  \item [(1)]
 If $D_1,D_2\in |H|$
are distinct then $D_1\cap D_2$ is a reduced
irreducible surface.
\item [(2)]
If $D_1,D_2,D_3\in |H|$ are linearly independent
the subscheme $D_1\cap D_2\cap D_3$ has pure
dimension $1$ and the Poincar\'e dual of the
fundamental cycle $[D_1\cap D_2\cap D_3]$ is
equal to $h^3$.
\end{itemize}
\end{prp}
\begin{proof}
(1): Assume that $\G\in Z^2(X)$ is an effective
non-zero algebraic cycle of pure codimension $2$.
Assume that
$$ cl(\G)=(sh^2+t(2q^{\vee}/5))\in H^4(X;\ZZ)/Tors,$$
where $cl(\G)$ is the image of the Poicar\'e dual
of the homology class represented by $\G$, and
$h:=c_1(H)$. Let $\s\in\G(\Omega^2_X)$ be a
symplectic form. Then
\begin{align*}
0< & \la cl(\G), h^2\ra=
\langle sh^2+t(2q^{\vee}/5),h^2\rangle=12s+20t,\\
0\le & \la cl(\G), (\s+\ov{\s})^2\ra=
\langle
sh^2+t(2q^{\vee}/5),(\s+\ov{\s})^2\rangle=(2s+10t)(\s+\ov{\s},\s+\ov{\s}).
\end{align*}
Since $(\s+\ov{\s},\s+\ov{\s})>0$ we get that
\begin{equation}\label{poscond}
3s+5t>0,\qquad s+5t\ge 0.
\end{equation}
Now let $D_1,D_2\in|H|$ be distinct. By Item~(2)
of Proposition~(\ref{hodgeprop}) we know that
$D_1\cap D_2$
 is a
subscheme of $X$ of pure codimension $2$
representing $h^2$. Assume that $D_1\cap D_2$ is
not reduced and irreducible: then we have an
equality of cycles $[D_1\cap D_2]=A+B$ with $A,B$
effective non-zero. By Item~(5) of
Proposition~(\ref{hodgeprop}) we have
$$cl(A)=xh^2+y(2q^{\vee}/5),\qquad
cl(B)=(1-x)h^2-y(2q^{\vee}/5)$$
with $2x,2y\in\ZZ$. Applying~(\ref{poscond}) we
get that
$$0<3x+5y<3,\qquad 0\le x+5y\le 1.$$
\lq\lq Eliminating $x$\rq\rq we get that
$$-3/5<2y<3/5.$$
Since $2y\in\ZZ$ we get that $y=0$ and hence
$cl(A)=xh^2$ with $0<x<1$. This contradicts
Item~(4) of Proposition~(\ref{hodgeprop}) and
proves Item~(1). Item~(2) follows immediately
from Item~(1).
\end{proof}
Let $B$ be the base-scheme of $|H|$, i.e.
\begin{equation}\label{luogobase}
 B:=\bigcap\limits_{D\in |H|}D.
\end{equation}
Item~(2) of Proposition~(\ref{intprop}) gives
that
\begin{equation}\label{tripleint}
  \dim B\le 1.
\end{equation}
We claim that
\begin{equation}\label{explicitrr}
h^0(\cO_X(nH))=\frac{1}{2}n^4+\frac{5}{2}n^2 +3,\quad n\in\NN_{+}.
\end{equation}
In fact applying H.-R.-R. and keeping in mind that all odd Chern classes of
$X$ vanish we get that for any $n\in\ZZ$ we have
\begin{equation}
\chi(\cO_X(nH))= \frac{1}{24}\left(\int_X
h^4\right) n^4+ \frac{1}{24}\left(\int_X
c_2(X)h^2\right) n^2 +\chi(\cO_X).
\end{equation}
By using~(\ref{c2formula}) and~(\ref{qdualint}) we get that
\begin{equation}\label{eulchar}
\chi(\cO_X(nH))=
\frac{1}{2}n^4+\frac{5}{2}n^2 +3,\quad n\in\ZZ.
\end{equation}
Since $K_X\cong\cO_X$ Kodaira vanishing gives that for $n>0$ we have
$h^0(\cO_X(nH))=\chi(\cO_X(nH))$. Thus~(\ref{explicitrr}) follows
from~(\ref{eulchar}). In particular we have $\chi(\cO_X(H))=6$.  We choose
once and for all an isomorphism
\begin{equation}\label{isomproj}
  |H|^{\vee}\overset{\sim}{\lra}\PP^5
\end{equation}
and we let
\begin{equation}\label{fmap}
    f\cl X\cdots>\PP^5
\end{equation}
be the rational map given by the composition
$X\cdots>|H|^{\vee}\overset{\sim}{\lra}\PP^5$.
Let
\begin{equation}\label{dilatoix}
\wt{X}:= Bl_B(X),
 \quad E\in Div(\wt{X})
\end{equation}
be the blow-up of the scheme $B$ and the
corresponding exceptional divisor respectively.
Let
\begin{equation}\label{effetilde}
  \wt{f}\cl\wt{X}\to \PP^5
\end{equation}
be the regular map which resolves the
indeterminacies of $f$. Let $Y:=Im(\wt{f})$; thus
$Y\subset\PP^5$ is closed and we have (abusing
notation) a dominant map
\begin{equation}\label{daxay}
  f\cl X\cdots> Y.
\end{equation}
We let $\deg f$ be the degree of the map above.
Let $Y_0$ be the interior of $\wt{f}(X\sm B)$ (we
may view $(X\sm B)$ as an open subset of
$\wt{X}$); thus $Y_0\subset Y$ is open and dense.
Let $X_0:=(X\sm B)\cap\wt{f}^{-1}(Y_0)$; thus
$X_0\subset X$ is open and dense. The restriction
of $\wt{f}$ to $X_0$ defines a regular surjective
map
\begin{equation}\label{effezero}
  f_0\cl X_0\to Y_0.
\end{equation}
\begin{prp}\label{nonredux}
Keep notation as above.
 Let $L\subset\PP^5$ be a  linear
subspace of codimension at most $2$. Then $L\cap
Y_0$ is reduced and irreducible and, if
non-empty, it has pure codimension equal to
$\cod(L,\PP^5)$.
\end{prp}
\begin{proof}
If $L=\PP^5$ there is nothing to prove. Assume
that $\cod(L,\PP^5)=1$. Let $D\in|H|$ be the
divisor corresponding to $L$
via~(\ref{isomproj}). Then $D\cap X_0=f_0^{*}L$;
since $X_0$ is open dense in $X$ and $f_0$ is
surjective the result follows from Item~(2) of
Proposition~(\ref{hodgeprop}). Assume that
$\cod(L,\PP^5)=2$ and write $L=L_1\cap L_2$ where
$L_1,L_2\subset\PP^5$ are hyperplanes. Let
$D_1,D_2\in|H|$ be the divisors corresponding to
$L_1,L_2$ via~(\ref{isomproj}). Then $D_1\cap D_2
\cap X_0=f_0^{*}L$; since $X_0$ is open dense in
$X$ and $f_0$ is surjective the result follows
from Item~(1) of Proposition~(\ref{intprop}).
\end{proof}
The following result is the first step towards
the proof that the manifold $X$ satisfies~(1)
or~(2) of Theorem~(\ref{mainthm1}).
\begin{prp}\label{varicasi}
Keep notation as above. One of the following
holds:
\begin{itemize}
\item[(1)] $\dim Y=3$ and $3\le \deg Y\le 6$. If
  $\dim Y=3$ and $\deg Y=6$ then $B$ is
  $0$-dimensional.
\item [(2)]
$\dim Y=4$, $\deg Y=2$.
\item [(3)]
$\dim Y=4$, $\deg Y=3$ and $\deg f= 3$.
\item [(4)]
$\dim Y=4$, $\deg Y=3$, $\deg f=4$ and $B=\es$.
\item [(5)]
$\dim Y=4$, $\deg Y=4$, $\deg f= 3$ and $B=\es$.
\item [(6)]
There exists a regular anti-symplectic involution $\phi\cl X\to X$ such that
$Y\cong X/\la\phi\ra$ and the quotient map $X\to X/\la\phi\ra$ is identified
with $f\cl X\to Y$. The $(\pm 1)$-eigenspaces of $H^2(\phi)$ are $\CC h$ and
$h^{\bot}$ respectively. The fixed locus of $\phi$ is a smooth irreducible
Lagrangian surface $F$ such that
\begin{equation}\label{chernancora}
  c_2(F)=192,\quad
  \cO_F(2K_F)\cong\cO_F(6H),\quad c_1(F)^2=360.
\end{equation}
\item [(7)]
$\dim Y=4$, $f\cl X\cdots>Y$ is birational and
$6\le \deg Y\le 12$.
\end{itemize}
\end{prp}
The rest of this section is devoted to the proof of the above proposition. We
let $d:=\deg Y$.
\begin{clm}\label{dimalmeno3}
Keeping notation as above, we have $\dim Y\ge 3$.
\end{clm}
\begin{proof}
This is a straightforward consequence of
Proposition~(\ref{nonredux}).  Suppose that $\dim
Y=1$. Since $Y$ is an irreducible non-degenerate
curve in $\PP^5$ we have $d\ge 5$. Let
$L\subset\PP^5$ be a generic hyperplane; since
$Y_0$ is open dense in $Y$ the intersection
$Y_0\cap L$ consists of $d$ points, contradicting
Proposition~(\ref{nonredux}).
 Now suppose that
$\dim Y=2$; since $Y$ is an irreducible
non-degenerate surface in $\PP^5$ we have $d\ge
4$. Let $L\subset \PP^5$ be a generic linear
subspace of codimension $2$; since $Y_0$ is open
dense in $Y$ the intersection $Y_0\cap L$
consists of $d$ points, contradicting
Proposition~(\ref{nonredux}).
\end{proof}
\n
 {\bf The case $\dim Y=3$.\/} We will show that~(1) holds.
 Since $Y$ is an irreducible non-degenerate $3$-fold
in $\PP^5$ we have $3\le\deg Y$. Let's prove that
$\deg Y\le 6$. Let $L,L',L''\subset \PP^5$ be
generic linearly independent hyperplanes. Then
the intersection $Y\cap L\cap L'\cap L''$ is
transverse and it consists of $d$ points
$p_1,\ldots,p_d\in Y_0$. Let $D,D',D''\in |H|$
correspond to $L,L',L''$ via~(\ref{isomproj}). By
Item~(2) of Proposition~(\ref{intprop}) the
scheme $D\cap D'\cap D''$ has pure dimension $1$.
Let $\G_{0,i}:=f_0^{-1}(p_i)$ and $\G_i$ be its
closure in $X$. We have
\begin{equation}\label{roberts}
  [D\cap D'\cap D'']=\G_1+\cdots+\G_d+\Si
\end{equation}
where $\Si$ is an effective $1$-cycle with  $supp
\Si\subset supp B$. (See~(\ref{fundcyc}) for the
notation $[D\cap D'\cap D'']$ .) Of course
$B\not=\es$ because $\dim Y<\dim X$, and
\begin{equation}\label{ovvio}
  \mbox{$\dim B=0$ if and only if $\Si=0$.}
\end{equation}
By~(3) of Proposition~(\ref{hodgeprop})
\begin{equation}\label{gammaclass}
  cl(\G_i)=m_i h^3/6,\quad m_i\in\NN_{+}.
\end{equation}
By Item~(2) of Proposition~(\ref{intprop}) the
$1$-cycle $[D\cap D'\cap D'']$ represents $h^3$
and hence~(\ref{roberts}) gives that
\begin{equation}\label{sordi}
  12=\la h,\G_1+\cdots+\G_d+\Si\ra=
 2\sum\limits_{i=1}^{d}m_i+\la h,\Si\ra\ge
 2d+\la h,\Si\ra.
\end{equation}
Since $h$ is ample and $\Si$ is effective we get
that $d\le 6$. Furthermore if $d=6$ then $\la
h,\Si\ra=0$ and hence $\Si=0$; by~(\ref{ovvio})
we get that $\dim B=0$.
 \vskip 3mm

\n
 {\bf The case $\dim Y=4$: elementary considerations.\/}
 Let $D,D',D'',D'''\in
|H|$ be linearly independent divisors. We will make
 some
 elementary considerations on the relation
 between the intersection number
 $\int_X h^4$ and the intersection
 $D\cap\cdots\cap D'''$. These facts will also be
 useful later on.  Let $L,L',L'',L'''\subset
\PP^5$ be the hyperplanes corresponding to
$D,D',D'',D'''$ via~(\ref{isomproj}). Let
$\wt{f}$ and $E$ be as in~(\ref{effetilde})
and~(\ref{dilatoix}) respectively; we can and
will assume that
\begin{equation}\label{evitae}
\dim(L\cap\cdots\cap L'''\cap Y)=0,\quad
 L\cap\cdots\cap L'''\cap\wt{f}(supp E)=\es.
\end{equation}
By Item~(2) of Proposition~(\ref{intprop}) the
intersection $D'\cap D''\cap D'''$ is of pure
dimension $1$.  There is a unique decomposition
\begin{equation}\label{simpledec}
  [D'\cap D''\cap D''']=\G+\Si
\end{equation}
with $\G,\Si$  effective $1$-cycles and
\begin{equation}\label{condizioni}
   \dim(supp (\G)\cap supp (B))\le 0, \quad
   supp\Si\subset supp B.
\end{equation}
From~(\ref{simpledec}) we get that
\begin{equation}\label{passetto}
 12=\int_X h^4=\deg(H\cdot(\G+\Si))=
 \deg(D\cdot\G)+\int_{\Si}h.
\end{equation}
By~(\ref{evitae}) the divisor $D$ intersects $\G$ in $d\cdot\deg f$ points
(counting multiplicities) outside $supp B$ and hence we have
\begin{equation}\label{sommetta}
 \deg(D\cdot\G)= d\cdot\deg f+\sum_{p\in supp B} mult_p(D\cdot\G).
\end{equation}
(The sum on the right is finite because of~(\ref{condizioni}).)
\begin{lmm}\label{productbound}
Keep notation as above. Assume that $\dim Y=4$. Then
\begin{equation}\label{tuttalpiu}
  \deg Y\cdot\deg f\le 12
\end{equation}
with equality if and only if $B=\es$.
\end{lmm}
\begin{proof}
Since $H$ is ample the integral appearing
in~(\ref{passetto}) is non-negative and
hence~(\ref{tuttalpiu}) follows
from~(\ref{passetto}) and~(\ref{sommetta}). It is
clear that if $B=\es$ then~(\ref{tuttalpiu})  is
an equality, we must prove the converse. Assume
that~(\ref{tuttalpiu}) is an equality.
By~(\ref{passetto}) and~(\ref{sommetta}) we have
$\Si=0$ and hence Equation~(\ref{simpledec})
gives that $supp\G\supset supp B$. Since $supp
D\supset supp B$ every $p\in supp B$ is contained
in $D\cap\G$. By~(\ref{sommetta}) we get that
$B=\es$.
\end{proof}

\n
 {\bf The case $\dim Y=4$ and $\deg f=1$.\/}  We must show
that~(7) holds. From Lemma~(\ref{productbound})
we get that $d\le 12$. One gets  the lower bound
$6\le d$ by adjunction. Explicitly, let
$\wt{Y}\subset\wt{\PP}^5$ be an embedded
resolution of $Y\subset\PP^5$: then
\begin{equation}\label{topform}
  h^0(K_{\wt{Y}})=1
\end{equation}
because $\wt{Y}$ is birational to $X$. On the
other hand by adjunction and vanishing of the
Hodge numbers
$h^{5,1}(\PP^5),h^{5,0}(\PP^5),h^{4,0}(\PP^5)$ we
get an isomorphism
\begin{equation}\label{adjform}
  H^0(K_{\wt{Y}})=H^0(I_Z(d-6)),
\end{equation}
where $Z\subset\PP^5$ is a subscheme supported on
$sing Y$. From~(\ref{topform}) we get that $6\le
d$. We have proved that if $\deg f=1$ then~(7)
holds.

\vskip 3mm

\n
 {\bf The case $\dim Y=4$ and $\deg f=2$.\/}
 Since $f\cl X\cdots>Y$ is generically a double
cover it defines a birational involution $\phi\cl
X\cdots>X$. We claim that $\phi$ is regular:
since $K_X\sim 0$ there exist closed subsets
$I_1,I_2\subset X$ of codimension at least $2$
such that $\phi$ restricts to a regular map
$(X\sm I_1)\to (X\sm I_2)$ and since
$H^{1,1}_{\ZZ}(X)=\ZZ h$ we have $\phi^{*}H\sim
H$; it follows by a well-known argument (see
\cite{huy}) that $\phi$ is regular. The map $f\cl
X\cdots> Y$ factors as
\begin{equation}\label{effedec}
  X\overset{\rho}{\lra} X/\la\phi\ra
  \overset{\ov{f}}{\cdots>}Y
\end{equation}
where $\rho$ is the quotient map. Since $\deg
f=2$ we have $\deg\ov{f}=1$, i.e.~$\ov{f}$ is
birational. We claim that
\begin{equation}\label{ynice}
  d=6,\quad\mbox{$\ov{f}$ is regular},
  \quad \dim (sing Y)\le 2.
\end{equation}
Let $\s$ be a symplectic form on $X$: since
$H^0(\Omega^2_X)=\CC\s$ and since $\phi$ is an
involution we have $\phi^{*}\s=\pm\s$ and hence
$\phi^{*}(\s\wedge\s)=\s\wedge\s$. Thus if $W$ is
any desingularization of $X/\la\phi\ra$ we have
$H^0(K_W)\not=0$. Since $\ov{f}$ is birational we
get that $H^0(K_{\wt{Y}})\not=0$ for any
desingularization $\wt{Y}\to Y$.
By~(\ref{adjform}) we get that $d\ge 6$, and
hence Lemma~(\ref{productbound}) gives that $d=6$
and that $B=\es$. Since $B=\es$ the map $\ov{f}$
is regular. Since $d=6$ we get that $\dim(sing
Y)\le 2$ - if $\dim(sing Y)=3$ then $sing Y$
certainly \lq\lq imposes conditions on
adjoints\rq\rq. We have proved~(\ref{ynice}).
Let's show that $\ov{f}$ is an isomorphism. The
fibers of $\ov{f}$ are finite because
$\ov{f}^{*}\cO_Y(1)$ is ample, $Y$ is normal
because it is a hypersurface smooth in
codimension $1$: this implies that the birational
map  $\ov{f}$ is an isomorphism. Let
$H^2_{\pm}(X)\subset H^2(X)$ be the
 $(\pm 1)$-eigenspace of $H^2(\phi)$
 respectively.  Then
 $h^2_{+}(X)$ is equal to $h^2(Y)$, which
is $1$ by Lefschetz' Hyperplane Section Theorem: since
$h$ belongs to
$H^2(\phi)_{+}$ we get that
\begin{equation}\label{accaduepiu}
  H^2(\phi)_{+}=\CC h.
\end{equation}
Since $\phi$ preserves Beauville's form $(,)$ we get that
\begin{equation}\label{accaduemeno}
 H^2(\phi)_{-}=h^{\bot}.
\end{equation}
In particular $\phi$ is anti-symplectic. Let's
prove that the fixed locus $F$ has the stated
properties. Since $F$ is the fixed locus of an
involution on a smooth manifold it is smooth.
Since $\phi$ is anti-symplectic $F$ has pure
dimension equal to $\dim X/2=2$, and $F$ is
Lagrangian. Let's prove that $F$ is irreducible.
Let $F=\bigcup\limits_{i\in I} F_i$ be the
decomposition into irreducible components. For
$i\in I$ let $cl(F_i)\in H^{2,2}_{\QQ}(X)$ be the
Poincar\'e dual of $F_i$; we claim that
\begin{equation}\label{classedieffe}
 cl(F_i)=k_i(15 h^2-c_2(X)),\quad k_i\in\QQ_{+}.
\end{equation}
In fact since $F_i$ is effective and Lagrangian we have
\begin{equation}\label{posenul}
 \int_X(cl(F_i)\wedge h^2)>0,\quad \int_X(cl(F_i)\wedge\s\wedge\ov{\s})=0.
\end{equation}
By Item~(6) of Proposition~(\ref{hodgeprop}) and by~(\ref{c2formula}) we have
\begin{equation}\label{espressione}
  cl(F_i)=(x_i h^2+y_i c_2(X)),\quad x_i,y_i\in\QQ.
\end{equation}
Substituting the above expression for $cl(F_i)$
in~(\ref{posenul}) and
applying~(\ref{collegamento})-(\ref{donpol})
and~(\ref{qdualint}) we get~(\ref{classedieffe}).
Now suppose that there exist two distinct
irreducible components $F_i,F_j$ of $F$. Then
$F_i\cap F_j=\es$ because $F$ is smooth and hence
by~(\ref{classedieffe}) we get that
\begin{equation}\label{prodzero}
 0=\int_X(cl(F_i)\wedge cl(F_i))=k_i k_j \int_X(15 h^2-c_2(X))^2.
\end{equation}
Thus $\int_X(15 h^2-c_2(X))^2=0$. On the other hand
using~(\ref{collegamento})-(\ref{donpol}) and~(\ref{qdualint}) we get that
\begin{equation}\label{belnum}
  \int_X(15 h^2-c_2(X))^2=1728,
\end{equation}
contradiction. This shows that $F$ is irreducible.
Let's prove that
\begin{equation}\label{cardif}
  c_2(F)=192.
\end{equation}
First we compute the Euler characteristic of $Y$. We have $b_i(Y)=\dim
H^i(\phi)_{+}$. Thus $b_i(Y)=0$ for odd $i$ and $b_2(Y)=1$
by~(\ref{accaduepiu}). By~(\ref{symmeq})
and~(\ref{accaduepiu})-(\ref{accaduemeno}) we get that $H^4(\phi)_{+}=\CC
(h\wedge h)\op Sym^2(h^{\bot})$ and hence $b_4(Y)=254$. Thus
\begin{equation}\label{cardiy}
  \chi(Y)=258.
\end{equation}
On the other hand the decompositions $X=(X\sm
F)\coprod F$ and $Y=(Y\sm \rho(F))\coprod
\rho(F)$
 give that
\begin{equation}\label{cardix}
  324=\chi(X)=2\chi(Y\sm \rho(F))+\chi(F).
\end{equation}
By~(\ref{cardiy}) we have $258=(\chi(Y\sm
\rho(F))+\chi(F))$; together with~(\ref{cardix})
this gives $\chi(F)=192$, i.e.~(\ref{cardif}).
Before proving the stated properties of $K_F$ we
show that
\begin{equation}\label{classevera}
cl(F)=5 h^2-\frac{1}{3} c_2(X).
\end{equation}
We have
\begin{equation}\label{normlagr}
\int cl(F)\wedge cl(F)=\int_F c_2(N_{F/X})=
\int_F c_2(\Omega^1_F)=192,
\end{equation}
where the second equality holds because $F$ is
Lagrangian and the third equality is given
by~(\ref{cardif}); replacing $cl(F)$ by the
right-hand side of~(\ref{classedieffe}) and
using~(\ref{belnum}) one gets~(\ref{classevera}).
Now let's prove that
\begin{equation}\label{formulacan}
  \cO_F(2K_F)\cong\cO_F(6H).
\end{equation}
Let $F':=\rho(F)$; thus $\rho\cl F\to F'$ is an
isomorphism. The embedding of $Y\cong
(X/\la\phi\ra)$ into $\PP^5$ defines by pull-back
an isomorphism
\begin{equation}\label{alquadro}
  \rho^{*}N_{F'/\PP^5}^{\vee}\cong
  Sym^2(N_{F/X}^{\vee}).
\end{equation}
Since $F$ is Lagrangian in $X$ we have
$N_{F/X}^{\vee}\cong\Theta_F$; substituting
in~(\ref{alquadro}) and taking determinants we
get an isomorphism
\begin{equation}\label{giudice}
  \rho^{*}\det(N_{F'/\PP^5})\cong\cO_F(3K_F).
\end{equation}
On the other hand the normal sequence for the
embedding $F'\hra\PP^5$ gives
\begin{equation}\label{piemme}
  \det(N_{F'/\PP^5})\cong\cO_{F'}(6)\ot\cO_{F'}(K_{F'}).
\end{equation}
Since $\rho$ is an isomorphism and
$\rho^{*}\cO_{F'}(1)\cong \cO_{F}(H)$ we get that
\begin{equation}\label{avvocato}
  \rho^{*}\det(N_{F'/\PP^5})\cong
  \cO_{F}(6H)\ot\cO_{F}(K_{F}).
\end{equation}
The above isomorphism together
with~(\ref{giudice}) gives~(\ref{formulacan}).
Finally to get $c_1(F)^2=360$
use~(\ref{formulacan}) and~(\ref{classevera})
together with~(\ref{collegamento})-(\ref{donpol})
and~(\ref{qdualint}). This finishes the proof
that if $\deg f=2$ then (6) holds.

We remark that we have the following stability
result for the $X$ satisfying~(6) of
Proposition~(\ref{varicasi}).
\begin{prp}\label{stabile}
Let $X$ be a numerical $(K3)^{[2]}$ and suppose
that there exist an anti-symplectic involution
$\phi\colon X\to X$ with quotient map $f\colon
X\to Y$ and an embedding $j\colon Y\hra\PP^5$
with $j(Y)$ a sextic hypersurface. Let
$H\sim\cO_Y(1)$. Let $X'$ be a small deformation
of $X$ for which $H$ remains of type $(1,1)$.
There is an involution $\phi'\colon X'\to X'$
which is a deformation of $\phi$ and letting
$f'\colon X'\to Y'$ be the quotient map there is
an embedding $Y'\hra\PP^5$ which deforms
$Y\hra\PP^5$. Furthermore $(f')^{*}\cO_{Y'}(1)$
is the divisor-class deformation of $H$.
\end{prp}
\begin{proof}
 Let $h:=c_1(H)$. Since $j(Y)$ is a sextic and
 $\deg f=2$ we have $\int_X h^4=12$. By
 Remark~(\ref{collegamento}) and
 Equation~(\ref{donpol}) we get that $(h,h)=2$.
 The invariant subspace $H^2(X)_{+}\subset H^2(X)$
 for the action of $H^2(\phi)$ contains $h$ and
 has rank $1$ because $H^2(Y)$ has rank $1$; thus
$H^2(X)_{+}=\CC h$. It follows that
$H^2(\phi)=R_h$ the reflection in the span of
$h$. The result then follows from
Proposition~(3.3) of~\cite{oginv}. (Notice that
in that proposition we have $0\in\cV$.)
\end{proof}
 \vskip 3mm

\n
 {\bf The case $\dim Y=4$ and $\deg f\ge 3$.\/}
By Lemma~(\ref{productbound})  we get
that one of (2), (3), (4), (5) holds.
 \vskip 3mm

\n
 We have proved Proposition~(\ref{varicasi}).
\section{Proof of
Theorem~(\ref{mainthm1})}\label{dimomainthm}
 \setcounter{equation}{0}
It suffices to prove that (1)-(5) of
Proposition~(\ref{varicasi}) cannot hold. We
assume that $f\cl X\cdots> Y$ satisfies on
of~(1), (2), ...~(5) of
Proposition~(\ref{varicasi}) and we reach a
contradiction. If $f\cl X\cdots> Y$ satisfies one
of~(1),...~(4) we show that either there exists a
linear subspace $L\subset\PP^5$ of codimension
$2$ such that $L\cap Y_0$ is not reduced and
irreducible  of pure codimension $2$,
contradicting Proposition~(\ref{nonredux}), or
the pull-back $f^{*}\cl H^4(Y)\to H^4(X)$ gives a
rational Hodge substructure of $H^4(X)$ which
does not exist by Proposition~(\ref{hodgeprop}).
If $f\cl X\cdots> Y$ satisfies~(5) of
Proposition~(\ref{varicasi}) and $\dim (sing
Y)=3$ then the first argument given above works.
If $f\cl X\cdots> Y$ satisfies~(5) of
Proposition~(\ref{varicasi}) and $\dim (sing
Y)\le 2$ then we show that the ramification
divisor of $f$ is the pull-back of a divisor on
$X$; since the ramification divisor is non-empty
this is absurd.
\subsection{(1) of
Proposition~(\ref{varicasi}) does not hold}
 \label{casouno}
 \setcounter{equation}{0}
We will prove the following result.
\begin{prp}\label{spezza}
Let $Y\subset\PP^5$ be an irreducible
non-degenerate linearly normal $3$-dimensional
subvariety of degree at most $6$.
\begin{itemize}
  \item [(1)]
 If $\deg Y\le 5$ then given an arbitrary non-empty
  subset $U\subset Y$ there exists a $3$-dimensional
  linear subspace $L\subset\PP^5$ such that $L\cap U$
  is reducible.
 \item [(2)]
If $\deg Y=6$ then there exists a $3$-dimensional
  linear subspace $L\subset\PP^5$ such that $L\cap Y$
  is not reduced or not irreducible.
\end{itemize}
\end{prp}
Granting the above proposition let's show that
(1) of Proposition~(\ref{varicasi}) does not
hold. The proof is by contradiction. First assume
that (1) of Proposition~(\ref{varicasi}) holds
with $\deg Y\le 5$. Clearly $Y$ is irreducible
non-degenerate and linearly normal and hence
Proposition~(\ref{spezza}) applies with $U:=Y_0$;
thus there exists a $3$-dimensional linear
subspace $L\subset\PP^5$ such that $L\cap Y_0$ is
reducible; this contradicts
Proposition~(\ref{nonredux}). This proves that we
cannot have $\deg Y\le 5$.
\begin{clm}\label{misterti}
Suppose that (1) of Proposition~(\ref{varicasi})
holds with $\deg Y=6$.  Then $Y_0=Y$.
\end{clm}
\begin{proof}
By Item~(1) of Proposition~(\ref{varicasi}) we
know that $\dim B=0$. Let $n$ be such that $nH$
is very ample and let $D\in |nH|$ be generic; in
particular since $\dim B=0$ we have  $D\subset
(X\sm B)=X_0$. It suffices to show that
\begin{equation}\label{tutto}
  f_0(D)=Y.
\end{equation}
Since $\dim Y=3$ the generic fiber of $f_0\cl
X_0\to Y$ is
 $1$-dimensional and hence its intersection with
$D$ consists of a finite set of points. Thus
$f_0(D)$ is $3$-dimensional. Since $f_0(D)$ is
closed in $Y$ and $Y$ is irreducible of dimension
$3$ we get~(\ref{tutto}).
\end{proof}
Now assume that (1) of
Proposition~(\ref{varicasi}) holds with $\deg
Y=6$; we will get to a contradiction. By
Claim~(\ref{tutto}) we have $Y_0=Y$. Since $Y$ is
irreducible non-degenerate and linearly normal
Proposition~(\ref{spezza})  applies and we get
that there exists a $3$-dimensional linear
subspace $L\subset\PP^5$ such that $L\cap Y_0$ is
not reduced or not irreducible, contradicting
Proposition~(\ref{nonredux}).
\vskip 3mm
 \n
 {\bf Proof of Proposition~(\ref{spezza}).}
 First consider the case of $Y$ a cone;
 thus $Y=J(p,\ov{Y})$ where $\ov{Y}$ is a surface
 with $\dim(span \ov{Y})=4$. (See~(\ref{joindef})
 for the notation $J(\cdot,\cdot)$.) Let
 $\ov{L}\subset(span \ov{Y})$ be a generic linear subspace of
 dimension $2$. Then $L:=J(p,\ov{L})$ is a
 $3$-dimensional linear subspace of $\PP^5$ and
\begin{equation}\label{joinjoin}
  L\cap Y= J(p,\ov{L}\cap\ov{Y}).
\end{equation}
Thus $L\cap U$ has $\deg\ov{Y}$ irreducible
components - they are open dense subsets of lines
through $p$. Since $\deg\ov{Y}=\deg Y\ge 3$ we
get that $L\cap Y$ is reducible. This proves the
proposition for $Y$ a cone. Now assume $Y$ is not
a cone. We prove Item~(1). Assume first that $Y$
is singular. Let $p\in sing(Y)$ and let $m$ be
its multiplicity. Let $A\subset\PP^5$ be a
hyperplane not containing $p$ and let
\begin{equation}\label{proietto}
  \rho\cl (Y\sm \{p\})\to A
\end{equation}
be projection from $p$. Let $Z:=Im(\rho)$ and let
$\ov{Z}$ be its closure. Since $Y$ is not a cone
$\ov{Z}$ is a hypersurface with $\deg
\ov{Z}=(\deg Y-m)$. Thus $\ov{Z}$ is a
hypersurface in $A\cong\PP^4$ of degree at most
$3$ and hence it is covered by lines. The image
$\rho(U\sm\{p\})\subset\ov{Z}$ contains an open
dense $V\subset Z$. Let $\ell\subset \ov{Z}$ be a
generic line: then $\ell\cap V$ is dense in
$\ell$. Let $q\in(V\sm\ell)$ be generic and let
$\ov{L}:=J(q,\ell)$. Thus $\ov{L}\subset A$ is a
plane and
\begin{equation}\label{elleci}
  \ov{L}\cap V=(\ell\cap V)\cup C
\end{equation}
where $C$ is an open dense subset of a line or of
a conic. (Notice that $\ov{L}\not\subset\ov{Z}$
because $\ell$ and $q$ are generic in $\ov{Z}$.)
Let $L:=J(p,\ov{L})$; this is a
 $3$-dimensional linear subspace of $\PP^5$.
We have
\begin{equation}\label{evidente}
  L\cap(\rho^{-1}V)=\rho^{-1}(\ov{L}\cap V)
\end{equation}
and hence $L\cap(\rho^{-1}V)$ is reducible
because of~(\ref{elleci}). Since $\rho^{-1}V$ is
an open subset of $U$ we get that $L\cap U$ is
reducible. Finally assume that $Y$ is smooth with
$\deg Y\le 5$. All smooth non-degenerate linearly
normal $3$-folds in $Y\subset\PP^5$ of degree at
most $5$ have been classified, see~\cite{ion}:
$Y$ is the Segre $3$-fold i.e.~$\PP^1\tm\PP^2$
embedded by
$\cO_{\PP^1}(1)\boxtimes\cO_{\PP^2}(1)$, or a
complete intersection of two quadric
hypersurfaces, or a quadric fibration, i.e.~it
fibers over $\PP^1$ with fibers which are
embedded quadric surfaces. In each case $Y$ is
covered by lines; it follows immediately that
Item~(1) of Proposition~(\ref{spezza}) holds for
$Y$. Now we prove Item~(2).   First assume that
$\dim(sing Y)=2$. Let $V\subset sing Y$ be a
$2$-dimensional component. We claim that
\begin{equation}\label{alpiuquattro}
  \deg V\le 4.
\end{equation}
In fact let $\Si\subset\PP^5$ be a generic
$3$-dimensional linear subspace: then
\begin{equation}\label{contiene}
  sing(\Si\cap Y)=\Si\cap sing Y\supset \Si\cap V
\end{equation}
and $|\Si\cap V|=\deg V$. Now $\Si\cap Y$ is an
irreducible non-degenerate curve in $\Si$, and
hence it has at most $4$ singular points.
Thus~(\ref{alpiuquattro}) follows
from~(\ref{contiene}). A straightforward argument
shows that any surface $V$ of degree at most $4$
contains a plane curve. Explicitely: If
$\dim(span V)= 2$ there is nothing to prove. If
$\dim(span V)\le 3$ intersect $V$ with a plane
contained in $span(V)$. If $\dim(span V)\ge 4$
and $V$ is singular the projection of $V$ from
$q\in(sing V)$ is a quadric surface $Q$; if
$\ell\subset Q$ is a line the intersection
$J(q,\ell)\cap V$ has dimension $1$. If
$\dim(span V)\ge 4$ and $V$ is smooth then
(see~\cite{ion}) $V$ is a rational scroll, a
complete intersection of quadric hypersurfaces in
a hyperplane of $\PP^5$ or the Veronese surface.
In the first two cases $V$ contains lines, in the
third case it contains conics. Thus we verified
that $V$ contains a plane curve $C$. Let
$L\subset\PP^5$ be the generic $3$-dimensional
linear space containing $C$: then $L\cap Y$ is a
reducible curve. This proves that Item~(2) holds
if $\dim(sing Y)=2$. Now assume that $\dim(sing
Y)=1$. Let $W\subset (sing Y)$ be a
$1$-dimensional component. If $\dim(span W)\le 2$
then $Y$ contains a plane curve and we are done.
Assume that $\dim(span W)\ge 4$. Then $\dim
((span W)\cap Y)\ge 2$ and hence there exists
$p\in ((span W)\cap (Y\sm W))$. Since curves are
never defective (see~\cite{ciro}) there exists a
$3$-secant plane of $W$ containing $p$, call it
$\Omega$. We claim that $\dim(\Omega\cap Y)\ge
1$. In fact if this is not the case then
$\dim(\Omega\cap Y)=0$ and hence the multiplicity
of the intersection $\Omega\cap Y$ is equal to
$\deg Y=6$: but the points in $\Omega\cap W$ give
a contribution of at least $6$ because $\Omega$
is $3$-secant to $W$ and $W\subset(sing Y)$, and
we have a contribution of at least $1$ from $p$,
for a total of at least $7$, contradiction. Thus
$Y$ contains a plane curve and we are done. We
are left with the case $\dim(span W)=3$. If $\dim
((span W)\cap Y)=2$ then $Y$ contains plane
curves and we are done. If $\dim ((span W)\cap
Y)=1$ let $L:=span W$; since $Y$ is singular
along $W$ the intersection $L\cap Y$ is  not
reduced along $W$. We have proved that Item~(2)
holds if $\dim(sing Y)\ge 1$. Now assume that
$\dim(sing Y)\le 0$. Let $\Lambda\subset\PP^5$ be
a generic hyperplane; thus $S:=\Lambda\cap Y$ is
a smooth non-degenerate (in $\Lambda$!) surface
of degree $6$. Since $\deg(S)\not=4$ we know that
$S$ is linearly normal (Severi) and we may apply
the known classification of such surfaces
(see~\cite{ion}): $S$ is the complete
intersection of a quadric and a cubic or it is a
Bordiga surface i.e.~the blow up of $\PP^2$ at
$10$ points embedded by the linear system of
plane quartics through the $10$ points. If $S$ is
a Bordiga surface it contains lines; if
$\ell\subset S$ is a line and $L\subset\PP^5$ a
generic $3$-dimensional linear subspace
containing $\ell$ the intersection $L\cap Y$ is
reducible. If $S$ is the complete intersection of
a quadric and a cubic then since $Y$ is linearly
normal the quadric hypersurface in $\Lambda$
containing $S$ lifts to a quadric hypersurface
$Q\subset\PP^5$ containing $Y$. There exist
$3$-dimensional linear spaces $L\subset\PP^5$
such that $L\cap Q$ is the union of $2$ planes;
if $L$ is a generic such space then $L\cap Y$ is
reducible. This finishes the proof of
Proposition~(\ref{spezza}).
\subsubsection{Comments}
One may ask the following: does there exist a
numerical $(K3)^{[2]}$ with an ample $H$ with
$(c_1(H),c_1(H))=2$ and $Y:=Im(f\cl X\cdots>|H|)$
of dimension strictly smaller than $4$? We do not
know of any such example however we do have
examples with $H$ big and nef such that $\dim
Y<\dim X$. (The case of big and nef divisors will
be needed in order to construct complete moduli
spaces.) An explicit example is the following.
Let $\pi\cl S\to\PP^2$ be a double cover ramified
over a smooth sextic; thus $S$ is a $K3$ surface.
Let $H_S:=\pi^{*}\cO_{\PP^2}(1)$ and let
$X:=M(0,H_S,0)$ be the Moduli space of
$H_S$-semistable rank-$0$ pure sheaves $G$ on $S$
with $c_1(G)=c_1(H_S)$ and $\chi(G)=0$: a typical
$G$ is given by $\iota_{*}\xi$ where $\iota\cl
C\hra S$ is the inclusion of a curve $C\in |H_S|$
and $\xi$ is a degree-$1$ line-bundle on $C$. It
is known that $X$ is a deformation of
$(K3)^{[2]}$ - see~\cite{yoshi}. There is a
Lagrangian fibration $\rho\cl X\to |H_S|$ mapping
$[G]\in M(0,H_S,0)$ to its support; the fiber
over $C\in |H_S|$ is $Jac^1(C)$ (suitably defined
if $C$ is singular). Thus on $X$ we have the
divisor class $F:=\rho^{*}\cO_{|H_S|}(1)$. We
also have a unique effective divisor $A$ on $X$
whose restriction to any Lagrangian fiber
$\rho^{-1}([C])\cong Jac^1(C)$ is the canonical
$\T$-divisor. Let $H:=A+2F$; a straightforward
argument shows that $(c_1(H),c_1(H))=2$ -
use~(\ref{hilbcase}). One can also show that $H$
is nef; since $\int_X c_1(H)^4=12$ we get that
$H$ is big. The image $Y=Im(f\cl X\cdots>|H|)$ is
the Veronese surface in $\PP^5$.
\subsection{(2) of
Proposition~(\ref{varicasi}) does not hold}
 \setcounter{equation}{0}
We assume that $Y\subset\PP^5$ is an irreducible
quadric hypersurface and we will get to a
contradiction. Since $Y_0\subset Y$ is open dense
in a quadric $4$-fold there exists a
$3$-dimensional linear subspace $L\subset\PP^5$
such that $L\cap Y_0$ is reducible; this
contradicts Proposition~(\ref{nonredux}).
\subsubsection{Comments}
There exist examples $(X,H)$ with $X$ a
deformation of $(K3)^{[2]}$ and $H$ an ample
divisors with $(c_1(H),c_1(H))=2$ such that
$Y=Im(f\cl X\cdots>|H|)$ is a quadric
hypersurface - see~(4.1) of~\cite{oginv}.
\subsection{(3) of
Proposition~(\ref{varicasi}) does not hold}
 \setcounter{equation}{0}
We will use the following elementary result.
\begin{prp}\label{tantipiani}
Let $Y\subset\PP^5$ be a cubic hypersurface
containing a $3$-dimensional linear space
$\Omega$. There exists a hyperplane
$Z\subset\PP^5$ containing $\Omega$ such that
$Z\cap Y$  is swept out by planes, i.e.~either
$Z\subset Y$ or $Z\cdot Y=\Omega+Q$ where
$Q\subset Z$ is a singular quadric hypersurface.
\end{prp}
\begin{proof}
Let $I\subset \Gr (3,\PP^5)\tm
|\cO_{\PP^5}(3)|\tm (\PP^5)^{\vee}$ be the set of
triples $(\Omega,Y,Z)$ where $\Omega\subset Y$
and $\Omega\subset Z$, let $J\subset\Gr
(3,\PP^5)\tm |\cO_{\PP^5}(3)|$ be the set of
couples $(\Omega,Y)$ where $\Omega\subset Y$ and
let
\begin{equation}\label{dimentica}
\begin{matrix}
  I & \overset{\rho}{\lra} & J \\
  (\Omega,Y,Z) & \mapsto & (\Omega,Y)
\end{matrix}
\end{equation}
be the forgetful map. Let $I^0\subset I$ be the
subset of triples $(\Omega,Y,Z)$ such that
$Z\cdot Y=\Omega+Q$ with $Q\subset Z$ a smooth
quadric hypersurface. We must show that
$\rho(I\sm I^0)=J$. The map $\rho$ is proper with
$1$-dimensional fibers, $J$ is irreducible and
$(I\sm I^0)$ is closed of codimension at most $1$
at every point; thus  it suffices to exhibit one
couple $(\Omega,Y)\in J$ such that
\begin{equation}\label{epimappa}
\rho^{-1}(\Omega,Y)\cap(I\sm I^0)\not=
\es,\rho^{-1}(\Omega,Y).
\end{equation}
Let $[X_0,\ldots,X_5]$ be homogeneous coordinates
on $\PP^5$. Let $\Omega=V(X_4,X_5)$ and
$Y=V(F\cdot X_4+G\cdot X_5)$ where
$F,G\in\CC[X_0,\ldots,X_5]$ are homogeneous of
degree $2$ with $F(X_0,\ldots,X_4,0)$ and
$G(X_0,\ldots,X_3,0,X_5)$ quadratic forms of rank
$4$ and $5$ respectively. Then $(\Omega,Y)$ is a
couple satisfying~(\ref{epimappa}).
\end{proof}
Now suppose that~(3) of Proposition~(\ref{varicasi}) holds i.e.~that $f\cl
X\cdots> Y$ is a map of degree $3$ and that $Y$ is a cubic hypersurface;  we
will arrive   at a contradiction. First we notice the following corollary of
Proposition~(\ref{tantipiani}).
\begin{crl}\label{nospaces}
Suppose that $f\cl X\cdots> Y$ is a map of degree
$3$ and that $Y$ is a cubic hypersurface. Then
$Y$ does not contain a $3$-dimensional linear
space.
\end{crl}
\begin{proof}
By contradiction. Let $\Omega\subset Y$ be a
$3$-dimensional linear space. By
Proposition~(\ref{tantipiani}) there exists a
hyperplane $Z\subset\PP^5$ containing $\Omega$
such that $Z\cap Y$ is swept out by planes. We
claim that $Z\cap Y_0\not=\es$. In fact let
$\wt{f}$ and $E$ be as in~(\ref{effetilde}) and
(\ref{dilatoix}) respectively; if $Z\cap Y_0=\es$
then $supp(\wt{f}^{*}Z)\subset
 supp(E)$, absurd.  Let $y\in Z\cap Y_0$;
 by Proposition~(\ref{tantipiani}) there exists a plane
$\Lambda\subset (Z\cap Y)$ with $y\in\Lambda$.
Now let $y'\in (Y_0\sm Z)$ and let
$L\subset\PP^5$ be the $3$-dimensional linear
space $L:=J(y',\Lambda)$. Then
\begin{itemize}
  \item[(a)] either $L\subset Y$, or
  \item [(b)] $L\cap Y=\Lambda\cup\G$ where $\dim\G=2$
  and $\G\ni y'$.
\end{itemize}
If Item~(a) holds then $L\cap Y_0$ is non-empty
$3$-dimensional  (notice: we do not
  know whether our \lq\lq original\rq\rq $3$-dimensional
  linear space $\Omega\subset Y$ intersects $Y_0$) and
  if Item~(b) holds
then $L\cap Y_0$ is reducible. In either case we
contradict Proposition~(\ref{nonredux}).
\end{proof}
Let $B$ be the base-scheme of $|H|$. We know
by~(\ref{tripleint}) that $\dim B\le 1$: we
consider separately the two cases $\dim B=0$ and
$\dim B=1$. Suppose first that $\dim B=0$. Then
the $1$-cycle $\Si$ appearing
in~(\ref{simpledec}) is zero and
hence~(\ref{passetto})-(\ref{sommetta}) give that
\begin{equation}\label{tuttiquanti}
  \sum_{p\in supp B} mult_p(D\cap D'\cap D''\cap D''')=3.
\end{equation}
This implies that $B$ is the disjoint union of
$0$-dimensional schemes $B_i$ each of which is
curvilinear (contained in a smooth curve) and
supported on a single point. Let $\ell_i$ be the
length of $B_i$; a straightforward computation
shows that $E=\sum_{i}\ell_i E_i$ with $E_i$ a
prime divisor such that $\pi(E_i)=supp B_i$.
(Recall that $\pi\cl\wt{X}\to X$ is the blow-up
of $B$.) Furthermore each $E_i$ is isomorphic to
$\PP^3$ and
$\wt{f}^{*}\cO_Y(1)\cong\cO_{E_i}(1)$. Thus
$\wt{f}(E_i)\subset Y$ is a $3$-dimensional
linear space; this is absurd by
Corollary~(\ref{nospaces}). Thus we are left with
the case $\dim B=1$.
\begin{prp}\label{nopoints}
Keep notation as above and assume that $\dim
B=1$. Then $B$ is a reduced, irreducible, local
complete intersection of pure dimension $1$.
Furthermore the following hold:
\begin{itemize}
  \item [(a)]
 Let $\Si$ be the $1$-cycle  appearing
in~(\ref{simpledec}). Then
  $\Si=[B]$, hence we may identify $\Si$ with
  $B$.
  \item [(b)]
 Let $\G$ be the $1$-cycle  appearing
in~(\ref{simpledec}). Then
  $supp(\G)$ intersects $\Si$ in a single
  point $p$, $supp(\G)$ is
  smooth at $p$ with tangent direction not
  contained in $\T_p \Si$,
  and the unique component
  of $supp(\G)$
  through $p$ appears with multiplicity 1 in in the cycle
  $\G$.
  \item [(c)]
  As $D,\ldots,D'''$ vary among divisors such
  that~(\ref{evitae}) holds the
  point of intersection $supp(\G)\cap \Si$ varies in
  $\Si$, i.e.~it
   is not constant.
\end{itemize}
\end{prp}
\begin{proof}
By~(\ref{passetto}) and~(\ref{sommetta}) we get
that
\begin{equation}\label{contributi}
\sum_{p\in supp B} mult_p(D\cap \G) + \int_{\Si}h
=3.
\end{equation}
Using  Item~(3) of Proposition~(\ref{hodgeprop})
we get that $cl(\Si)=m h^3/6$ for some
 positive integer $m$; in fact $m$ is  strictly
positive or else $\dim B=0$.
 Thus $\int_{\Si}h=2m$ and by~(\ref{contributi})
we get that
\begin{equation}\label{classesigma}
  \int_{\Si}h=2,\quad cl(\Si)=h^3/6.
\end{equation}
Furthermore, again by Item~(3) of
Proposition~(\ref{hodgeprop}), we get that
$supp(\Si)$ is irreducible and that the
multiplicity of $\Si$ equals $1$, i.e.~$\Si$ is
equal to the reduced irreducible curve
$supp(\Si)$. Since $H$ is ample the scheme
$D'\cap D''\cap D'''$ is connected and hence
$supp(\G)\cap \Si\not=\es$; let $p\in
supp(\G)\cap \Si$. Since $\Si\subset B$ we have
$p\in supp B$; thus~(\ref{contributi}) gives that
\begin{equation}\label{moltuno}
  mult_p(D\cdot \G)=1,\quad
  supp(B)\cap supp(\G)=\{p\}.
\end{equation}
This proves Item~(b) because $\Si\subset
supp(B)$. Furthermore
 since $supp(B)\subset supp(\G)\cup\Si$ we get that
 $supp(B)=supp(\Si)$. On the other hand $B$ is a
subscheme of $D\cap\cdots\cap D'''$ and away from
$supp(\G)$ the latter scheme coincides with the
reduced and irreducible l.c.i.~$\Si$; thus $B$ is
a reduced and irreducible l.c.i.~away from $p$.
This shows that if Item~(c) holds then also
Item~(a) holds. We prove Item~(c) arguing by
contradiction. If Item~(c) is false then
$supp(\G)\cap \Si=\{p\}$ for a fixed $p\in\Si$
whenever~(\ref{evitae}) holds. Let
$\Lambda\subset|H|$ be the $\PP^3$ spanned by
$D,\ldots,D'''$ and let
\begin{equation}\label{lambdapi}
  \Lambda_p:=\{ Z\in\Lambda|\ mult_p Z\ge 2\}.
\end{equation}
Since all $Z\in\Lambda$ contain $\Si$ we get that
$\Lambda_p$ is a linear subspace of $\Lambda$
with $\cod(\Lambda_p,\Lambda)\le 3$ and hence
$\Lambda_p$ is not empty because $\dim\Lambda=3$.
Renaming $D,\ldots,D'''$ we may assume that
$D\in\Lambda_p$. Since $p\in supp(\G)\cap
supp(B)$ Equation~(\ref{contributi}) gives that
\begin{equation}\label{contraddizione}
  4\le \sum_{p\in supp B} mult_p(D\cap \G) +
  \int_{\Si}h=3,
\end{equation}
absurd.
\end{proof}
Let's show that if $\dim B=1$ then we get a
contradiction. By Proposition~(\ref{nopoints})
the exceptional divisor is a $\PP^2$-fibration
$E\to B$ and $\wt{f}$ embeds each fiber
$\pi^{-1}(p)$ over $p\in B$ as a plane in
$\PP^5$. We claim that $\wt{f}(E)$ is a
$3$-dimensional linear subspace of $\PP^5$. In
fact let $L',L'',L'''\subset\PP^5$ be the
hyperplanes corresponding to $D',D'',D'''$
via~(\ref{isomproj}). By
Proposition~(\ref{nopoints}) the divisors
$\wt{f}^{*}(L'),\wt{f}^{-1}(L''),
 \wt{f}^{*}(L'''), E$ intersect transversely in a
 single point. Thus $\dim\wt{f}(E)=3$, because if
 we had $\dim\wt{f}(E)=2$ then the intersection would be
 either empty or of dimension $1$. Furthermore since
 $L',L'',L'''$ are generic hyperplanes we
get that $\wt{f}(E)$ has degree $1$, i.e.~it is
$3$-dimensional  linear space. Since
$\wt{f}(E)\subset Y$ this contradicts
Proposition~(\ref{nospaces}). This completes the
proof that Item~(3) of
Proposition~(\ref{varicasi}) does not hold.
\subsection{(4) of
Proposition~(\ref{varicasi}) does not hold}
 \setcounter{equation}{0}
We will prove that  our map $f\cl X\cdots> \PP^5$
cannot be a degree-$4$ regular map onto an
irreducible cubic hypersurface  $Y\subset\PP^5$.
The proof is by contradiction. We assume that we
have $f\cl X\to Y$ a finite regular map of degree
$4$ onto a cubic $4$-fold $Y\subset\PP^5$ and we
reach a contradiction. If $Y$ is smooth a
straightforward argument shows that $f^{*}H^4(Y)$
is a non-existant Hodge substructure of $H^4(X)$
- see Subsubsection~(\ref{casoliscio}). The proof
that $Y$ cannot be a singular cubic $4$-fold is
more involved: it will follow from some results
on singular cubic $4$-folds which should be of
independent interest. Let $Y\subset\PP^5$ be an
arbitrary singular cubic hypersurface: for $p\in
sing(Y)$ we let
\begin{equation}\label{essepidef}
 S_p:=\{\ell\in\Gr(1,\PP^5)|\ p\in\ell\subset Y\}.
\end{equation}
The definition above is set-theoretic but of
course $S_p$ has a natural structure as subscheme
of $\Gr(1,\PP^5)$. We will prove the following
result.
\begin{prp}\label{autaut}
Let $Y\subset\PP^5$ be a singular cubic
hypersurface and let $p\in sing Y$. Then
\begin{itemize}
  \item [(1)]
either $Y$ contains a plane or
  \item [(2)]
$Y$ has isolated quadratic singularities,   the
scheme $S_p$ is a reduced, normal surface not
containing lines nor conics, with Du Val
singularities\footnote{See Ch.4
of~\cite{kollmori} for definition and properties
of Du Val singularities.} The minimal
desingularization of $S_p$ is a $K3$ surface
$\wt{S}_p$.
\end{itemize}
\end{prp}
The proof of Proposition~(\ref{autaut}) goes as
follows. If $Y$ is reducible then $Y$ satisfies
Item~(1) trivially - and of course it does not
satisfy Item~(2). If $Y$ is a cone then $Y$
certainly does not satisfy Item~(2), and it
satisfies Item~(1) by the following elementary
result.
\begin{lmm}\label{facile}
Let $Y\subset\PP^5$ be a cubic hypersurface which
is a cone. Then $Y$ contains a plane.
\end{lmm}
\begin{proof}
We have $Y=J(p,\ov{Y})$ where $\ov{Y}$ is a cubic
hypersurface in $\PP^4$. Thus $\ov{Y}$ contains a
line $\ell$ and hence $Y$ contains the plane
$J(p,\ell)$.
\end{proof}
We are left with the case of $Y$ an irreducible singular cubic $4$-fold which is
not a cone, i.e.~every singular point of $Y$ is quadratic. In
Subsubsections~(\ref{tredim})-(\ref{duedim})-(\ref{unodim}) we will prove that
if $\dim (sing Y)$ is equal to $3$, $2$, $1$ respectively then $Y$ contains  a
plane; the basic (well-known) observation is that the line joining two distinct
singular points of $Y$ is contained in $Y$. Thus if $(\dim Y)\ge 1$ then
Item~(1) of Proposition~(\ref{autaut}) holds - and Item~(2) does not hold by
hypothesis. In Subsubsection~(\ref{zerodim}) we prove that if $Y$ has isolated
quadratic singularities then either Item~(1) or Item~(2) of
Proposition~(\ref{autaut}) holds. It is elementary that~(1) and~(2) cannot both
hold; the hard part is to show that if~(1) does not hold then  $S_p$ has Du Val
singularities - the remaining statements of~(2) are straightforward with the
exception of the assertion about the minimal desingularization of $S_p$, this
follows from the fact that the singularities are Du Val. First we prove by
explicit computation that the singularities of $S_p$ which are $\ul{\text{not}}$
lines joining $p$ to another singular point of $Y$ are Du Val. Then by analyzing
the relation between $S_p$ and $S_{p'}$ for $p'\not=p$ we are able to get that
$S_p$ is Du Val also at the points $span(p,p')$ for $p'\in sing Y$. This will
complete the proof of Proposition~(\ref{autaut}). In order to prove that (4) of
Proposition~(\ref{varicasi}) does not hold we will need a result on the (mixed)
Hodge structure of a cubic $4$-fold $Y$ satisfying Item~(2) of
Proposition~(\ref{autaut}). Let
\begin{equation}\label{filtpeso}
  \ldots \subset W_{3} H^4(Y)\subset W_{4} H^4(Y)=H^4(Y)
\end{equation}
be Deligne's weight filtration~\cite{deligne}. In
particular
\begin{equation}\label{doppiavutre}
  W_3 H^4(Y)=\ker(H^4(Y)\overset{H^4(\zeta)}{\lra}
   H^4(\wt{Y}))
\end{equation}
where $\zeta\cl \wt{Y}\to Y$ is any
desingularization, see Proposition~(8.5.2)
of~\cite{deligne}. Thus $W_3 H^4(Y)$ is in the
kernel of the intersection form on $H^4(Y)$ and
hence the intersection form is well-defined on
$Gr_4^W H^4(Y):=H^4(Y)/W_3 H^4(Y)$. Let $p\in
sing Y$; since we are assuming that $Y$ satisfies
Item~(2) of Proposition~(\ref{autaut}) we know
that $\wt{S}_p$ is a $K3$ surface. Let
$T(\wt{S}_p)\subset H^2(\wt{S}_p;\ZZ)$ be the
transcendental lattice of $\wt{S}_p$ i.e.
\begin{equation}\label{trasc}
  T(\wt{S}_p):=\{\alpha\in H^2(\wt{S}_p;\ZZ)|\
  \alpha\bot H^{1,1}_{\ZZ}(\wt{S}_p)\}.
\end{equation}
Then
\begin{equation}\label{compltrasc}
  T(\wt{S}_p)_{\CC}:=T(\wt{S}_p)\ot_{\ZZ}\CC\subset
  H^2(\wt{S}_p)
\end{equation}
is a sub-Hodge structure of level $2$ with
\begin{equation}\label{accaduez}
  h^{2,0}(T(\wt{S}_p)_{\CC})=
  h^{0,2}(T(\wt{S}_p)_{\CC})=1,\quad
 1\le h^{1,1}(T(\wt{S}_p)_{\CC})\le 19.
\end{equation}
The following result will be proved in
Subsubsection~(\ref{firenze}).
\begin{prp}\label{eccemappa}
Suppose that a cubic hypersurface $Y\subset\PP^5$
satisfies Item~(2) of Proposition~(\ref{autaut}).
Then there is a morphism of type $(1,1)$ of Hodge
structures
\begin{equation}\label{unouno}
  \g\cl T(\wt{S}_p)_\CC\lra Gr_4^W H^4(Y).
\end{equation}
If $\eta,\theta\in T(\wt{S}_p)_\CC$ then
\begin{equation}\label{chiave}
  \int_Y \g(\eta)\wedge\gamma(\theta)=
  -\int_{\wt{S}_p}\eta\wedge\theta.
\end{equation}
\end{prp}
Granting
Propositions~(\ref{autaut})-(\ref{eccemappa})
let's prove that it is impossible to have $f\cl
X\to Y$ finite of degree $4$ onto a singular
cubic $4$-fold. Assume that such an $f$ exists;
since $f$ is regular $Y=Y_0$ and hence $Y$ does
not contain planes by
Proposition~(\ref{nonredux}). Let $p\in sing Y$.
By Propositions~(\ref{autaut})-(\ref{eccemappa})
we have the morphism of type $(1,1)$ of Hodge
structures $\g$ of~(\ref{unouno}). Composing $\g$
with $f^{*}$ we get a morphism of type $(1,1)$ of
Hodge structures
\begin{equation}
T(\wt{S}_p)_\CC\overset{f^{*}\circ \g}{\lra}
H^4(X).
\end{equation}
Let $\eta,\theta\in T(\wt{S}_p)_{\CC}$;
by~(\ref{chiave}) we have
\begin{equation}\label{nomadi}
\int_X f^{*}\g(\eta)\wedge f^{*}\gamma(\theta)=
  -4\int_{\wt{S}_p}\eta\wedge\theta.
\end{equation}
Since the restriction to $T(\wt{S}_p)_\CC$ of the
intersection form on $H^2(\wt{S}_p)$ is
non-degenerate we get that $f^{*}\circ \g$ is
injective. Thus $Im(f^{*}\circ \g)$ is a rational
Hodge substructure of $H^4(X)$ with Hodge numbers
$h^{p,q}=h^{p-1,q-1}(T(\wt{S}_p)_{\CC})$.
By~(\ref{accaduez}) this contradicts Item~(4) of
Proposition~(\ref{hodgeprop}).

In the last subsusbsection we comment on the
possibility that $f\cl X\to Y$ is of degree $4$
onto a cubic when one drops one of the hypotheses
of Proposition~(\ref{varicasi}).
\subsubsection{(4) of Proposition~(\ref{varicasi})
with $Y$ smooth does
not hold}\label{casoliscio}
We assume that $f\cl X\to Y$ with $Y\subset\PP^5$
a smooth cubic hypersurface, $f$ finite of degree
$4$ and we get to a contradiction. Since $\deg
f=4$ we have
\begin{equation}\label{perquattro}
  \la f^{*}\alpha, f^{*}\beta\ra_X=4\la \alpha,\beta\ra_Y,
  \quad \alpha,\beta\in H^4(Y)
\end{equation}
where $\la,\ra_X$ and $\la,\ra_Y$ are the intersection forms on $H^4(X)$ and
$H^4(Y)$ respectively. Thus $f^{*}\cl H^4(Y)\to H^4(X)$ is an injection of
rational Hodge structures. Let
$$H^4(Y)_{prim}:=\{\alpha\in H^4(Y)|\ \alpha\wedge
c_1(\cO_Y(1))=0\}$$
be the primitive cohomology of $Y$: this a rational sub Hodge structure of
$H^4(Y)$.
 Since $\dim
H^4(Y)_{prim}=22$ Item~(4) of Proposition~(\ref{hodgeprop}) gives that
$f^{*}H^4(Y)_{prim}=\CC h\ot h^{\bot}$. Thus
\begin{equation}\label{uguali}
 f^{*}H^4(Y;\QQ)_{prim}=\QQ h\ot h_{\QQ}^{\bot}
\end{equation}
where $h_{\QQ}^{\bot}:=h^{\bot}\cap
H^{2}(X;\QQ)$. Let
$\cB=\{\alpha_1,\ldots,\alpha_{22}\}$ be a
$\ZZ$-basis of $H^4(Y;\ZZ)_{prim}$. Let
$\cQ_{\cB}$ be the matrix of the restriction of
$\la,\ra_Y$ to $H^4(Y;\ZZ)_{prim}$ in the basis
$\cB$. Since $\la,\ra_Y$ is unimodular and $\deg
Y=3$ we have
\begin{equation}\label{valoretre}
  |\det(\cQ_{\cB})|=3.
\end{equation}
Let
$\cB':=\{f^{*}\alpha_1,\ldots,f^{*}\alpha_{22}\}$;
by~(\ref{uguali}) we know that $\cB'$ is a
$\QQ$-basis of $\QQ h\ot h^{\bot}_{\QQ}$. Let
$\cQ_{\cB'}$ be the matrix of the restriction of
$\la,\ra_X$ to $\QQ h\ot h^{\bot}_{\QQ}$ in the
basis $\cB'$; by
~(\ref{valoretre})-(\ref{perquattro}) we have
\begin{equation}\label{disctreper}
  |\det(\cQ_{\cB'})|=3\cdot 2^{44}.
\end{equation}
Now let $\{\beta_1,\ldots,\beta_{22}\}$ be a
$\ZZ$-basis of
 $h^{\bot}_{\ZZ}:=
H^2(X;\ZZ)\cap h^{\bot}$; then
$\cB'':=\{h\beta_1,\ldots,h\beta_{22}\}$ is a
$\QQ$-basis of $\QQ h\ot h^{\bot}_{\QQ}$. Let
$\cQ_{\cB''}$ be the matrix of the restriction of
$\la,\ra_X$ to $\QQ h\ot h^{\bot}_{\QQ}$ in the
basis $\cB''$. By Remark~(\ref{collegamento}) one
gets (use also Lemma~(\ref{tuttiequiv})) that
\begin{equation}\label{ventiquattro}
 |\det(\cQ_{\cB''})|=2^{24}.
\end{equation}
Since both $\cB'$ and $\cB''$ are $\QQ$-bases of $\QQ h\ot h^{\bot}_{\QQ}$
the determinants appearing in
Equations~(\ref{disctreper})-(\ref{ventiquattro}) must represent the same
class in $\QQ^{*}/(\QQ^{*})^2$. This is visibly false, contradiction.
\subsubsection{$Y$ a singular cubic $4$-fold:
elementary considerations}
\label{preliminari}
Let $Y\subset\PP^5$ be a singular cubic
hypersurface. Suppose that $p,q\in sing Y$ are
distinct points: $span( p,q)$ and $Y$ intersect
with  multiplicity at least $2$ at $p$ and at $q$
hence by B\'ezout we get that $span( p,q)\subset
Y$. Thus for a subset $W\subset sing Y$  we have
\begin{equation}\label{utile}
 chord(W)\subset Y
\end{equation}
where $chord (W)\subset\PP^5$ is the subvariety
swept out by the chords of $W$ i.e.
\begin{equation}\label{cordev}
  chord (W):=\text{closure of $\{span(p,q)|\ p,q\in
   W,\ \ p\not=q\}$.}
\end{equation}
Now assume that $Y$ is irreducible and reduced, $p\in sing Y$ and $Y$ is not a
cone with vertex $p$. Choose homogeneous coordinates $[X_0,\ldots,X_4,Z]$ on
$\PP^5$ such that $p=[0,\ldots,0,1]$. We have
\begin{equation}\label{eqesplicita}
  Y=V(F(X_0,\ldots,X_4)Z+G(X_0,\ldots,X_4))
\end{equation}
where $F,G$ are homogeneous non-zero of degrees
$2$ and $3$ respectively. We have
\begin{equation}\label{conoquadrica}
  \PP(C_p Y)=V(F(X_0,\ldots,X_4))\subset
  \PP^4_{[X_0,\ldots,X_4]}=\PP(\Theta_p Y)=\PP(\Theta_p \PP^5).
\end{equation}
Let
\begin{equation}\label{psiconpi}
  \psi_p\cl Y\cdots> \PP(\Theta_p \PP^5)
\end{equation}
be projection from $p$. The map $\psi_p$ is
birational: letting $\ul{X}:=X_0,\ldots,X_4$ the
inverse of $\psi_p$ is given by
\begin{equation}\label{invdipsi}
\begin{matrix}
 \qquad \qquad &\PP(\Theta_p \PP^5)& \overset{\psi_p^{-1}}{\cdots>}& Y \\
 & & & \\
 \qquad \qquad & [\ul{X}] & \mapsto &
 [F(\ul{X})X_0,\ldots,F(\ul{X})X_4,-G(\ul{X})]
\end{matrix}
\end{equation}
The indeterminacy locus  of $\psi_p^{-1}$ is
clearly the set of lines through $p$ contained in
$Y$ (see~(\ref{essepidef})). Using the
coordinates introduced above we see that the
natural inclusion $S_p\subset \PP(\Theta_p
\PP^5)$ is given by
\begin{equation}\label{essepi}
S_p=V(F,G)\subset
  \PP^4_{[\ul{X}]}=\PP(\Theta_p \PP^5).
\end{equation}
Notice that since $Y$ is irreducible, reduced and
not a cone with vertex $p$ the polynomials
$F,G\in\CC[\ul{X}]$ have no common factors and
hence
\begin{equation}\label{intercomp}
  \text{$S_p$ is a complete
  intersection of $\PP(C_p Y)$
   and a cubic hypersurface.}
\end{equation}
Formula~(\ref{invdipsi}) says that $\psi_p^{-1}$
is defined by the linear system $|I_{S_p}(3)|$ on
$\PP(\Theta_p \PP^5)$. Since $I_{S_p}(3)$ is
globally generated we get that the resolution of
indeterminacies of $\psi_p$ defines an
isomorphism
\begin{equation}\label{duescoppi}
  \wt{\psi}_p\cl Bl_p Y\overset{\sim}{\lra}
  Bl_{S_p}\PP(\Theta_p \PP^5).
\end{equation}
 We will need to relate properties of $Y$ and of
$S_p$. A first observation: if $y\in sing
(Y\sm\{p\})$ then $span( p,y)\subset Y$
by~(\ref{utile}) and hence
\begin{equation}\label{dentro}
  \psi_p(sing Y\sm\{p\})\subset S_p\subset\PP(C_p Y).
\end{equation}
\begin{prp}\label{singsing}
Suppose that $Y\subset\PP^5$ is a singular
reduced and irreducible cubic hypersurface, that
$p\in sing Y$ and that $Y$ is not a cone with
veretx $p$.
\begin{itemize}
  \item[(1)]
  If $y\in sing(Y\sm\{p\})$ then $s:=\psi_p(y)\in
  sing (S_p)$. If
  $span(p,y)\subset sing(Y)$ then
  $\dim\Theta_s(S_p)=4$, in particular $\PP(C_p Y)$ is
  singular at $s$. If $span(p,y)\not\subset sing(Y)$ then
  $\PP(C_p Y)$ is smooth at $s$.
 \item [(2)]
 Let $s\in sing (S_p)$ and assume that
 $\dim\Theta_s(S_p)=4$. Then $Y$ is singular at all points
 of the line corresponding to $s$.
\item [(3)]
Let $s\in sing (S_p)$ and assume that
 $\dim\Theta_s(S_p)=3$. If $\PP(C_p Y)$ is smooth at $s$
there exists a unique $y\in sing(Y\sm\{p\})$ such
that $\psi_p(y)=s$. If $\PP(C_p Y)$ is singular
at $s$ there is no $y\in sing(Y\sm\{p\})$ such
that $\psi_p(y)=s$.
 \item [(4)]
 $Y$ contains a plane if and only if $S_p$
 contains a line or a conic.
\end{itemize}
\end{prp}
\begin{proof}
Let $[X_0,\ldots,X_4,Z]$ be homogeneous coordinates on $\PP^5$ with
$p=[0,\ldots,0,1]$; thus we have~(\ref{eqesplicita})-(\ref{essepi}). Let
$y=[a_0,\ldots,a_4,b]\in\PP^5\sm\{p\}$: thus
\begin{equation}\label{eccoqui}
  \psi_p(y)=[a_0,\ldots,a_4]=[\ul{a}].
\end{equation}
Differentiating the defining equation of $Y$ we
get that $y\in sing(Y\sm\{p\})$ if and only if
\begin{equation}\label{jacob}
b\cdot\frac{\partial F}{\partial x_i}(\ul{a})+
 \frac{\partial G}{\partial x_i}(\ul{a})=0
 \quad i=0,\ldots,4,
 \quad\text{and}\quad  F(\ul{a})=0.
\end{equation}
(1): From the two equations above we get that
$G(\ul{a})=0$ (we already noticed this), and
hence the first equation shows that $s\in
sing(S_p)$. Assume that for a fixed
$\ul{a}\not=(0,\ldots,0)$ the first equation
holds with an arbitrary choice of $b$: then both
$V(F)$ and $V(G)$ are singular at $s$ and this
proves the second statement. Assume that for a
fixed $\ul{a}\not=(0,\ldots,0)$ the first
equation holds for some but not for all choices
of $b$: then $V(F)$ is smooth at $s$ and this
proves the third statement. Items~(2)-(3) are
proved by similar elementary considerations. Now
let's prove Item~(4). Assume that $Y$ contains a
plane $L$. If $p\in L$ then $\psi_p(L\sm\{p\})$
is a line contained in $S_p$. If $p\notin L$ then
$\Lambda:=\psi_p(L)$ is a plane in
$\PP^4_{[\ul{X}]}$. The restriction of
$\psi_p^{-1}$ to $\Lambda$ is the linear system
$|I_{\Lambda\cap S_p}(3)|$. Since
$\psi_p^{-1}(\Lambda)=L$ is a plane we get that
necessarily $\Lambda\cap S_p$ is a conic in
$\Lambda$; thus $S_p$ contains a conic. The proof
of the converse is similar.
\end{proof}
\subsubsection{Proof of Proposition~(\ref{autaut})
for $Y$ with $\dim(sing Y)=3$}\label{tredim}
As shown in the introduction to the subsection we
may assume that $Y$ is reduced, irreducible and
not a cone. Let $Y\subset\PP^5$ be a reduced and
irreducible cubic hypersurface with $\dim(sing
Y)=3$. The intersection of $Y$ and a generic
plane is a singular reduced and irreducible cubic
curve and hence it has exactly one singular
point. Thus $sing Y$ has exactly one
$3$-dimensional irreducible component, call it
$V$, and $V$ is a linear space. Thus $Y$ contains
(many) planes.
\subsubsection{Proof of Proposition~(\ref{autaut})
for $Y$ with $\dim(sing Y)=2$}\label{duedim}
 $Y$ is necessarily reduced and irreducible. We may also
 assume that $Y$ is
not a cone by Lemma~(\ref{facile}). Assume that there exists a
$2$-dimensional irreducible component $V$  of $sing Y$ with $\dim(span(V))\le
4$. Then $chord(V)=span
  (V)$ and hence by~(\ref{utile}) $Y$ contains a linear subspace of dimension at
  least $2$. Now assume that every
  $2$-dimensional irreducible component $V$
   of $sing Y$ is non-degenerate.
By~(\ref{utile}) we get that $\dim(chord(V))\le
4$, i.e.~the non-degenerate surface
$V\subset\PP^5$ is defective: a classical result
of Severi (see~\cite{ciro}) states that $V$ is
either a cone over a degree-$4$ rational normal
curve or the Veronese surface. One verifies
easily that in both cases $chord(V)$ is a cubic
hypersurface in $\PP^5$ and hence $Y=chord(V)$.
If $V$ is a cone over a degree-$4$ rational
normal curve then $chord(V)$ is itself a cone,
excluded by hypothesis. If $V$ is a Veronese
surface let $\psi\cl\PP^2\overset{\cong}{\to} V$
be an isomorphism with
$\psi^{*}\cO_{V}(1)\cong\cO_{\PP^2}(2)$; if
$\ell\subset\PP^2$ is a line then $\psi(\ell)$ is
a conic spanning a plane contained in $chord(V)$.
Thus $chord(V)=Y$ contains a plane.
\subsubsection{Proof of Proposition~(\ref{autaut})
for $Y$ with $\dim(sing Y)=1$}\label{unodim}
$Y$ is necessarily reduced and irreducible. We
may also
 assume that $Y$ is
not a cone by Lemma~(\ref{facile}). Let $(sing
Y)^1$ be the union of $1$-dimensional irreducible
components of $sing Y$.
 Choose $p\in (sing Y)^1$ such that
\begin{equation}\label{puntosing}
  \mbox{$(sing Y)^1$ is smooth at $p$.}
\end{equation}
Let $S_p$ be the set of lines in $Y$ through $p$
- see~(\ref{essepidef}). Assume first that $S_p$
is not reduced or that it is reducible.
By~(\ref{intercomp})  we get that there exists a
surface $T\subset S_p$ of degree at most $3$ and
hence $S_p$ contains  a line $\ell$. The lines in
$\PP^5$ parametrized by points of $\ell$ sweep
out a plane contained in $Y$, and we are done.
Now assume that $S_p$ is reduced and irreducible:
let's prove that
\begin{equation}\label{limitesup}
  \deg(sing Y)^1\le 5.
\end{equation}
Let $\psi_p\cl Y\cdots>\PP(\T_p \PP^5)$ be projection from $p$. Let $(sing
S_p)^1$ be the union of $1$-dimensional irreducible components of $sing(S_p)$ -
notice that $sing(S_p)$ has dimension at most $1$ becuse $S_p$ is a reduced
surface. By Proposition~(\ref{singsing}) the closure of
  $\psi_p(sing Y\sm\{p\})$ is an irreducible
  component of $(sing S_p)^1$ and hence
\begin{equation}\label{natavota}
 \deg\overline{\psi_p(sing Y\sm\{p\})}\le
 \deg(sing S_p)^1.
\end{equation}
 We claim that
\begin{equation}\label{alpiuqu}
  \deg(sing S_p)^1\le 4.
\end{equation}
In fact let $\Lambda\subset\PP(\Theta_p \PP^5)$
be a generic $3$-dimensional linear space; thus
$S_p\cap \Lambda$ is irreducible.
By~(\ref{essepi}) $S_p\cap \Lambda$ is a complete
intersection  of a quadric and a cubic in
$\Lambda\cong\PP^3$ and hence it has arithmetic
genus $4$; since it is irreducible we get that it
has at most $4$ singular points.
Inequality~(\ref{alpiuqu}) follows because
$sing(S_p\cap \Lambda)=sing(S_p)\cap\Lambda$. By
Assumption~(\ref{puntosing}) we have
\begin{equation}\label{menouno}
  \deg\psi_p((sing Y)^1)=\deg(sing Y)^1-1.
\end{equation}
Inequality~(\ref{limitesup}) follows
from~(\ref{menouno})-(\ref{natavota})-(\ref{alpiuqu}).
Thus if $S_p$ is reduced and irreducible one of
the following holds:
\begin{itemize}
  \item [(I)]
  $(sing Y)^1$ contains a line.
  \item [(II)]
 There is an irreducible component $\G$ of $(sing
 Y)^1$ with $2\le \dim(span(\G))\le 3$.
  \item [(III)]
 There is an irreducible component $\G$ of $(sing
 Y)^1$ with $\dim(span(\G)) = 4$ and $4\le\deg (\G)\le 5$.
  \item [(IV)]
 $(sing Y)^1$ is the rational normal curve of
 degree $5$ in $\PP^5$.
\end{itemize}
We will examine (I) through (IV) separately and
we will show in each case that $Y$ contains a
plane. (I): Let $\ell\subset sing Y$ be a line.
We will prove that there exists a plane
$\Lambda\subset Y$ containing $\ell$. Let
$[X_0,\ldots,X_5]$ be homogeneous coordinates on
$\PP^5$ such that $\ell=V(X_0,\ldots,X_3)$. Since
$Y$ is singular along $\ell$ we have
$$Y=V(A\cdot X_4+B\cdot X_5+C)$$
where $A,B,C\in\CC[X_0,\ldots,X_3]$ are homogeneous with $\deg A=\deg B=2$ and
$\deg C=3$. There exists a point
$$[a_0,\ldots,a_3]\in V(A,B,C)\subset\PP^3_{[X_0,\ldots,X_3]}.$$
The plane
$$\Lambda:=\{[\lambda a_0,\ldots,\lambda a_3,\mu,\theta]|
\ [\lambda,\mu,\theta]\in\PP^2\}$$
is contained in $Y$. (II): Since $\dim(span(\G))\le 3$ we have
$chord(\G)=span(\G)$. By~(\ref{utile}) we know that $Y\supset span(\G)$. Since
by hypothesis $\dim(span(\G))\ge 2$ we get that $Y$ contains a plane. (III):
First we prove the following.
\begin{lmm}\label{razquattro}
Let $Y\subset\PP^5$ be a reduced and irreducible
cubic hypersurface such that $sing Y$ contains an
irreducible curve $\G$  with $\dim(span( \G)) =
4$ and $4\le\deg (\G)\le 5$. Then $\G$ is a
degree-$4$ rational normal curve and $Y\cap(span(
\G))$ is the cubic $3$-fold $chord(\G)$.
\end{lmm}
\begin{proof}
By~(\ref{utile}) $chord(\G)\subset Y$. The
intersection $Y\cap(span(\G))$ is a
  hypersurface because $Y$ is reduced and irreducible. Since $chord(\G)$ is a
  hypersurface in $span(\G)$ we get that
\begin{equation}\label{nonminore}
  3=\deg(Y\cap span(\G))\ge\deg(chord(\G)),
\end{equation}
with equality only if $(Y\cap span(\G))=(chord(\G))$. From our hypotheses we get
that either $\G$ is a degree-$4$ rational normal curve in $span(\G)$ or it has
degree $5$ and arithmetic genus at most $1$. A straightforward computation shows
 that
 $$\deg(chord\G)=
  \begin{cases}
    3 & \text{if $\deg\G=4$}, \\
    6 & \text{if $\deg\G=5$ and $p_a(\G)=0$},\\
     5  & \text{if $\deg\G=5$ and $p_a(\G)=1$}.
  \end{cases}$$
The result follows from the above formulae and~(\ref{nonminore}).
\end{proof}
 Now fix a degree-$4$ rational normal curve $\G\subset\PP^5$. If it were true that
  $chord(\G)$ contains a plane we would be done; unfortunately this is not the
  case. Let $\cI_{\G}\subset\cO_{\PP^5}$ be the ideal sheaf of $\G$;
  thus $|\cI_{\G}^2(3)|$ is the linear system of cubic hypersurfaces
  $Y\subset\PP^5$ which are singular at each point of $\G$. Before formulating the
  next result we remark that the $3$-fold $chord(\G)$ contains lines which are
  not chords of $\G$.
\begin{prp}\label{pianoper}
Keep notation as above, and let $Y\in
|\cI_{\G}^2(3)|$. Then $Y$ contains a
$1$-dimensional family of planes $\Lambda$ such
that $\Lambda\cap span(\G)$ is a chord of $\G$.
\end{prp}
\begin{proof}
Let $Z\subset\G^{(2)}\tm\Gr(2,\PP^5)$ be the
subset defined by
$$Z:=\{(p+q,\Lambda)|\ \Lambda\supset\ov{p,q}\},$$
where $\ov{p,q}=span(p,q)$ if $p\not=q$ and $\ov{ p,p}=T_p \G$. Projecting  $Z$
to the first factor we get that $Z$ is smooth irreducible and
\begin{equation}\label{dimcinque}
  \dim Z=5.
\end{equation}
Let $(p+q,\Lambda)\in Z$: we let
$\cI_{p+q,\Lambda}\subset\cO_{\Lambda}$ be the
ideal sheaf of the subscheme $\{p,q\}$ (reduced
structure) if $p\not=q$ and of the length-$2$
subscheme supported at $p$ with tangent direction
$\Theta_p\G\subset\Theta_p\Lambda$ if $p=q$. Let
$F\to Z$ be a vector-bundle with fiber
$H^0(\cI_{p+q,\Lambda}(2))$ over $(p+q,\Lambda)$.
Of course $F$ is only defined modulo
tensorization by a line-bundle on $Z$: any choice
of $F$ is good for our argument. We have
\begin{equation}\label{rangoeffe}
  \rk F=4.
\end{equation}
Let  $Y=[P]\in |\cI_{\G}^2(3)|$ where
$P\in\CC[X_0,\ldots,X_5]$ is homogeneous of
degree $3$. Let $z=(p+q,\Lambda)\in Z$ and
$\tau_{z}\in H^0(\cO_{\Lambda}(1))$ be an
equation of the line $\ov{p,q}\subset\Lambda$. We
have
\begin{equation}\label{restralambda}
P|_{\Lambda}=\tau_z\ot\sigma_{z,P},\quad
\sigma_{z,P}\in H^0(I_{p+q,\Lambda}(2)).
\end{equation}
Letting $\pi\cl Z\tm |\cI_{\G}^2(3)|\to Z$ be the projection the above equation
gives that there is a section $\s\in H^0(\pi^{*}F\ot\cL)$, where  $\cL\to Z\tm
|\cI_{\G}^2(3)|$ is a suitable line-bundle, such that
\begin{equation}\label{sezmult}
  \s(z,[P])=c\cdot \s_{z,P},\quad c\in\CC^{*}.
\end{equation}
Let $W:=(\s)$ be the locus of zeroes of $\s$. Letting $\rho\cl Z\tm
|\cI_{\G}^2(3)|\to |\cI_{\G}^2(3)|$ be the projection we have
$$\rho(W)=\{Y\in |\cI_{\G}^2(3)|\  |\ \text{$\exists \Lambda\subset Y$ a plane with
$\Lambda\cap span(\G)$ a chord of $\G$.}\}$$
\begin{clm}\label{toro}
Keep notation as above. There exist $(z_0,Y_0)\in
W$ and an open $U\subset Z\tm |\cI_{\G}^2(3)|$
containing $(z_0,Y_0)$ such that $U\cap W\cap
\rho^{-1}(Y_0)$ is purely $1$-dimensional.
\end{clm}
\begin{proof}
As is easily checked there exists a smooth $Q\in
|\cI_{\G}(2)|$. Since $\G$ is cut out by quadrics
we may assume that
\begin{equation}\label{noncord}
  Q\not\supset chord(\G).
\end{equation}
let $Y_0:=Q+span(\G)$; clearly $Y_0\in
|\cI_{\G}^2(3)|$. Before choosing $z_0$ we notice
that
$$\Si_Q:=\{p+q\in\G^{(2)}|\ \ov{p,q}\subset Q\}$$
is $1$-dimensional because of~(\ref{noncord}).
Let $p_0+q_0\in\Si_Q$. There exist two planes
$\Lambda\subset Q$ which contain $\ov{p_0,q_0}$,
let $\Lambda_0$ be one of them: we set
$z_0:=(p_0+q_0,\Lambda_0)$. We let $U\subset Z\tm
|\cI_{\G}^2(3)|$ be the open subset given by
$$U:=\{(p+q,\Lambda,Y)|\ \Lambda\not\subset span(\G)\}.$$
 One easily checks that with these choices
 the claim holds.
\end{proof}
Let's finish the proof of the proposition.
By~(\ref{rangoeffe}) we get that $\cod(W,Z\tm
|\cI_{\G}^2(3)|)\le 4$, and thus
by~(\ref{dimcinque})
\begin{equation}\label{dimdoppiav}
\dim W\ge 1+\dim|\cI_{\G}^2(3)|.
\end{equation}
By Claim~(\ref{toro}) the fibers of $\rho$
restricted to $W\cap U$ have dimension at most
$1$ in a neighborhood of $Y_0$ and hence
$$\dim\rho(W)=\dim|\cI_{\G}^2(3)|.$$
Since $\rho$ is proper and $|\cI_{\G}^2(3)|$ is
irreducible we get that
$\rho(W)=|\cI_{\G}^2(3)|$, i.e.~every $Y\in
|\cI_{\G}^2(3)|$ contains a plane intersecting
$span(\G)$ in a chord of $\G$. Furthermore the
set of such planes has dimension at least $1$
because every fiber of $\rho|_W$ has dimension at
least $1$ by~(\ref{dimdoppiav}) and because every
plane in $\PP^5$ intersects $\G$ in a finite set
of points.
\end{proof}
\n
 (IV): Let $\G\subset\PP^5$ be a rational normal curve
 of degree $5$. We will
explicitly construct cubic hypersurfaces
$Y\subset\PP^5$ with $\G\subset sing(Y)$; by
construction these cubics are ruled by planes.
Then we will prove that every $Y\in
|\cI^2_{\G}(3)|$ is one of the cubics that we
constructed; thus every cubic satisfying (IV)
contains a plane - actually a $2$-dimensional
family.  Let $L\to\G$ be \lq\lq the\rq\rq
degree-$1$ line-bundle. Given a degree-$3$ linear
system $G$ of dimension $2$ on $\G$ i.e. $G\in
|L^{\ot 3}|^{\vee}$, we let
\begin{equation}\label{gidef}
  Y_G:=\bigcup\limits_{p_1+p_2+p_3\in G}
  \ov{p_1,p_2,p_3}
\end{equation}
be the variety swept out by the planes spanned by
divisors parametrized by $G$ - of course if
$p_1=p_2=p$ and $p_3\not=p$ then
$\ov{p_1,p_2,p_3}:=J(T_p \G,p_3)$
 and if $p_1=p_2=p_3=p$ then $\ov{p_1,p_2,p_3}$ is the the projective
 osculating plane to $\G$ at $p$. One easily checks that $Y_G$ is a hypersurface
 and that $sing(Y_G)=\G$.
 Furthermore $Y_G$ is a cone with vertex $p$ if and only if $p\in\G$ and $G=
 p+|L^{\ot 2}|$; if this is the case then $Y_G=\la p,chord(\G_p)\ra$
 where $\G_p\subset\PP(\Theta_p(\PP^5))$ is the projection of $\G$ from $p$. Since
 $\G_p$ is a degree-$4$ rational normal curve
 $chord(\G_p)$ is a cubic $3$-fold and hence we get that $\deg(Y_G)=3$ whenever
 $G$ has a base point. Since $\deg(Y_G)$ is independent of $G$ we get that $Y_G$
 is a cubic hypersurface for all $G\in|L^{\ot 3}|^{\vee}$.
 Thus we have defined an injection
\begin{equation}\label{birigata}
\begin{matrix}
|L^{\ot 3}|^{\vee} & \hra & |\cI^2_{\G}(3)| \\
G & \mapsto & Y_G
\end{matrix}
\end{equation}
\begin{prp}\label{eccecubica}
Keep notation as above. The map~(\ref{birigata})
is an isomorphism.
\end{prp}
\begin{proof}
Let $p\in\G$ and let
$\Si_p\subset|\cI^2_{\G}(3)|$ be the linear
subspace of cubics which are cones with vertex
$p$. Let $G_p:=(p+|L^{\ot 2}|)\in |L^{\ot
3}|^{\vee}$; a straightforward argument shows
that
\begin{equation}\label{dimensig}
  \Si_p=\{Y_{G_p}\}.
\end{equation}
Now let's prove that
\begin{equation}\label{codimtre}
  \cod(\Si_p,|\cI^2_{\G}(3)|)\le 3.
\end{equation}
Let $U\ni p$ be an open affine space containing
$p$; associating to $Y\in|\cI^2_{\G}(3)|$ an
affine cubic equation of $Y\cap U$ we may
identify $H^0(\cI^2_{\G}(3))$ with a
sub-vector-space $A\subset\CC[U]$. If
$Y\in|\cI^2_{\G}(3)|$ then $Y$ is singular at
$p$; thus $p$ is a critical point of $\phi$ for
all $\phi\in A$. Associating to $\phi\in A$ its
Hessian at $p$ we get a linear map
\begin{equation}\label{hessiano}
\begin{matrix}
A & \overset{\cH}{\lra} & Sym_2(\Omega^1_p(\PP^5)) \\
\phi & \mapsto & \text{Hessian of $\phi$ at $p$.}
\end{matrix}
\end{equation}
Since $\Si_p=\PP(\ker\cH)$ it suffices to prove
that
\begin{equation}\label{pochequad}
  \dim(Im\cH)\le 3.
\end{equation}
Let $Q\in\PP(Im\cH)$; we may view $Q$ as a
quadric hypersurface in $\PP^5$ with vertex at
$p$. Since cubics in $|\cI^2_{\G}(3)|$ are
singular at all points of $\G$ the quadric $Q$ is
singular at all points of $T_p \G$. Moreover $Q$
contains all the lines $\ov{p,q}$ for $q\in\G$
because such lines are contained in any
$Y\in|\cI^2_{\G}(3)|$. Hence projecting $Q$ from
the line $T_p \G$ we get a quadric
$\ov{Q}\subset\PP(N_{T_p \G,\PP^5})$ containing
the degree-$3$ rational normal curve $\ov{\G}$
obtained projecting $\G$ from $T_p \G$. The
linear system of quadrics in $\PP(N_{T_p
\G,\PP^5})\cong\PP^3$ containing $\ov{\G}$ has
(projective) dimension $2$ and hence we
get~(\ref{pochequad}). This
proves~(\ref{codimtre}).
 By~(\ref{dimensig}) we get
that $\dim|\cI^2_{\G}(3)|\le 3$. Since the map
of~(\ref{birigata}) is injective and since
$\dim|L^{\ot 3}|^{\vee}=3$ we get the
proposition.
\end{proof}
\subsubsection{Proof of Proposition~(\ref{autaut})
for $Y$ with $\dim (sing Y)=0$}\label{zerodim}
We assume that $Y\subset\PP^5$ is a singular
cubic hypersurface with isolated singularities
and that $p\in sing Y$. First let's show that
Items~(1) and~(2) of Proposition~(\ref{autaut})
are mutually exclusive. Assume that $Y$ contains
a plane. We may assume that $Y$ is reduced and
irreducible because otherwise $Y$ does not have
isolated singularities. If $Y$ is  a cone with
vertex $p$ then $\dim(S_p)=3$ and hence Item~(2)
of Proposition~(\ref{autaut}) does not hold. If
$Y$ is not a cone with vertex $p$ then by
Item~(4) of Proposition~(\ref{singsing}) we know
that $S_p$ contains a line or a conic and hence
Item~(2) of Proposition~(\ref{autaut}) does not
hold. This shows that Items~(1) and~(2) of
Proposition~(\ref{autaut}) can not both hold. It
remains to prove that if $Y$ does not contain
planes then Item~(2) of
Proposition~(\ref{autaut}) holds. By
Lemma~(\ref{facile}) we know that $Y$
 is not a cone, i.e.~it has quadratic singularities.
 By~(\ref{intercomp})
 we know that $S_p$ is a surface.
We say that a surface is {\it Du Val} if it is reduced,
normal with Du Val
singularities. We notice the following
\begin{equation}\label{bastaduval}
  (\text{$S_p$ is Du Val})\implies
  \text{Item~(2) of Proposition~(\ref{autaut}) holds.}
\end{equation}
In fact by~(\ref{intercomp}) we know that
$S_p$ is an intersection of a
quadric and a cubic in $\PP^4$  and hence by simultaneous resolution
  of Du Val singularities it follows that the minimal desingularization
  $\wt{S}_p$ is a deformation of a smooth intersection of a quadric and a
  cubic in $\PP^4$. Since a smooth intersection of a quadric and a
  cubic in $\PP^4$ is a $K3$ surface we get that $\wt{S}_p$ is a $K3$. Thus
  it remains to show that $S_p$ is Du Val.
\begin{clm}\label{pocopoco}
Let $Y\subset\PP^5$ be a singular cubic
hypersurface with isolated quadratic
singularities. Assume that $Y$ does not contain
any plane. Let $q\in sing(Y)$. Then:
\begin{itemize}
  \item [(1)]
   $\dim(sing \PP(C_q Y))\le 1$.
  \item [(2)]
$S_q\subset\PP(\Theta_q\PP^5)$ is a reduced and
irreducible normal complete intersection of
$\PP(C_q Y)$ and a cubic hypersurface. $S_q$ has
hypersurface singularities (embedding dimension
$3$).
\end{itemize}
\end{clm}
\begin{proof}
(1): Suppose that $\dim(sing \PP(C_q Y))\ge 2$.
Then $\PP(C_q Y))$ is the union of two
hyperplanes in $\PP(\Theta_q\PP^5)\cong\PP^4$ or
a double hyperplane, and hence $S_q$ is the union
of two cubic surfaces or a double cubic surface.
In either case $S_q$ contains a line,
contradicting Item~(4) of
Proposition~(\ref{singsing}). (2):  By Item~(4)
of Proposition~(\ref{singsing}) $S_q$ contains no
lines and hence by Item~(1) we get that $S_q\cap
sing \PP(C_q Y)$ is empty or finite. This fact
together with Items~(1)-(2)-(3) of
Proposition~(\ref{singsing}) gives that $S_q$ is
reduced normal and that the  embedding dimension
of $S_q$ is equal to $3$ at every singular point.
Furthermore~(\ref{intercomp}) gives that $S_q$ is
a complete intersection as stated. $S_q$ is
connected because it is a complete intersection:
since $S_q$ is normal we get that it is
irreducible.
\end{proof}
Keep notation as above. By the above claim $S_p$
is a reduced and normal surface with locally
trivial dualizing sheaf $\om_{S_p}$ (actually
$\om_{S_p}$ is globally trivial by adjunction).
It remains to prove that the singularities of
$S_p$ are Du Val, i.e.~that given any $s\in
sing(S_p)$ there exists a desingularization of
$s$, call it $\e_s\cl T_s\to S_p$, such that
$\om_{T_s}\cong\e_s^{*}\om_{S_p}$. Let
\begin{equation}\label{puntising}
  |sing(Y)|=k+1.
\end{equation}
Let $q\in sing(Y)$: we write
$sing(Y)=\{q,q_1,\ldots,q_k\}$. The line
$span(q,q_i)$ for $1\le i\le k$ is contained in
$Y$ and hence it is parametrized by a point
$[span(q,q_i)]\in S_q$.  Let
\begin{equation}\label{fuorising}
  U_q:=S_q\sm\{[span(q,q_1,)],\ldots,[span(q,q_k,)]\}.
\end{equation}
\begin{prp}\label{duralex}
Let $Y\subset\PP^5$ be a cubic with isolated
quadratic singularities and assume that $Y$ does
not contain any plane. Let $q\in sing(Y)$. The
singularities of $U_q$ are Du Val.
\end{prp}
\begin{proof}
By Proposition~(\ref{singsing}) we know that
$U_q$ is smooth away from $sing(\PP(C_q Y))\cap
S_q$. Thus we must prove that $S_q$ has a Du Val
singularity at all $s\in sing(\PP(C_q Y))\cap
S_q$. Choose such an $s$. By
Claim~(\ref{pocopoco}) we get that $S_q$ is the
complete intersection of $\PP(C_q Y)$ and a cubic
$\Xi\subset\PP(\Theta_q\PP^5)$ which is smooth at
$s$. From $\dim(sing\PP(C_q Y))\le 1$ one easily
gets that $mult_s(S_q)=2$. Let $\pi\cl
\wt{S}_q\to S_q$ be the blow-up of $s$. Since
$S_q$ has a hypersurface singularity of
multiplicity $2$ at $s$ we have
$\pi^{*}\om_{S_q}\cong\om_{\wt{S}_q}$. Thus it
suffices to prove that
\begin{equation}\label{bastaquesto}
  \text{$\wt{S}_q$ has Du Val singularities along
  $\pi^{-1}(s)$.}
\end{equation}
Since $\PP(C_q Y)$ is singular at $s$ it is the
join $J(s,Q)$ where $Q\subset\PP(\Theta_q \PP^5)$
is a quadric surface not containing $s$. By
Item~(1) of Claim~(\ref{pocopoco}) we know that
$Q$ is either smooth or the cone over a smooth
conic. By Item~(4) of
Proposition~(\ref{singsing}) we know that  $S_q$
contains no lines and hence projection from $s$
defines a regular finite map $\psi\cl\wt{S}_q\to
Q$ of degree $2$. We describe explicitly $\psi$.
Let  $\ul{X}:=X_0,\ldots,X_3$. Choose projective
coordinates $[\ul{X},Z]$ on $\PP(\Theta_q\PP^5)$
so that $s=[0,\ldots,0,1]$ and $span(Q)=V(Z)$;
thus $[\ul{X}]$ are projective coordinates on
$span(Q)$. We have
\begin{equation}\label{eccoequa}
\PP(C_q Y)=V(F),\qquad \Xi=V(A Z^2+B Z +C)
\end{equation}
where $F,A,B,C\in\CC[\ul{X}]$ are homogeneous of
degrees $2$, $1$, $2$ and $3$ respectively. Since
$S_q$ contains no lines we have
\begin{equation}\label{intvuota}
  \PP^3_{[\ul{X}]}\supset V(F,A,B,C)=\es.
\end{equation}
Since $\wt{S}_q$ is normal the branch divisor of $\psi\cl\wt{S}_q\to Q$ is
the reduced effective divisor $D(\psi)\in Div(Q)$ defined by
\begin{equation}\label{branchdiv}
  D(\psi)=V(F,B^2-4 A\cdot C)\subset
  Q=V(F)\subset\PP^3_{[\ul{X}]}.
\end{equation}
In general suppose that $V$ is a smooth surface, $W$ is a normal surface and
$\vf\cl W\to V$ is a double cover branched over the effective reduced divisor
$D\in Div(V)$. One can get a desingularization $\wh{W}$ of $W$ by
constructing an embedded resolution $\wh{D}$ of $D$ in a suitable blow-up
$\wh{V}$ of $V$ and taking a double cover $\wh{W}\to\wh{V}$ branched over
$\wh{D}$ and a suitable sum of components of the exceptional divisors: from
this construction one easily gets the following criterion.
\begin{crt}\label{duvalse}
Keep notation as above. Let $w\in W$ and $v:=\vf(w)$. Suppose that
 $mult_v (D)\le 3$ and moreover that if $mult_v (D)=3$ the strict transform of
$D$ in $Bl_v(V)$ intersects the exceptional divisor in at least two distinct
points. Then $W$ has a Du Val singularity at $w$.
\end{crt}
 Now let's prove~(\ref{bastaquesto}). Let $t\in\pi^{-1}(s)$
and let $[\ul{e}]=\psi(t)$. We have
\begin{equation}\label{keane}
  [\ul{e}]\in\psi(\pi^{-1}(s))=
  V(F,A)\subset Q\subset\PP^3_{[\ul{X}]}.
\end{equation}
If $B(\ul{e})\not=0$ then by~(\ref{keane})
and~(\ref{branchdiv}) we get that $[\ul{e}]\notin
D(\psi)$. Thus a neighborhood of $t$ in
$\wt{S}_q$ is isomorphic to a neighborhood of
$[\ul{e}]$ in $Q$. Since $Q$ has Du Val
singularities we get that $\wt{S}_q$ is Du Val at
$t$. Thus we may assume from now on  that
\begin{equation}\label{binullo}
B(\ul{e})=0.
\end{equation}
By~(\ref{intvuota}) we have
\begin{equation}\label{cinonzero}
  C(\ul{e})\not=0.
\end{equation}
We treat separately the two cases:
\begin{itemize}
  \item [(1)]
 $Q$ is smooth at $[\ul{e}]$.
  \item [(2)]
$Q$ is singular at $[\ul{e}]$.
\end{itemize}
(1): If $V(A)$ is transverse to $Q=V(F)$ at $[\ul{e}]$ then
by~(\ref{cinonzero}) we get that $D(\psi)$ is smooth at $[\ul{e}]$ and hence
$\wt{S}_q$ is Du Val at $t$ - actually smooth. If $V(A)$ is tangent to $Q$ at
$[\ul{e}]$ we distinguish two cases: $Q$ smooth and $Q$ singular. If $Q$ is
smooth then $V(A,F)$ is the union of two distinct lines through $[\ul{e}]$
and we get from~(\ref{cinonzero}) and~(\ref{branchdiv}) that $D(\psi)$ has a
quadratic singularity at $[\ul{e}]$: thus $\wt{S}_q$ is Du Val at $t$ by
Criterion~(\ref{duvalse}). If $Q$ is singular then $V(A,F)$ is a \lq\lq
double line\rq\rq supported on $\ell:=span([\ul{e}],sing Q)$. If $V(B)$ is
singular at $[\ul{e}]$ or if it is smooth at $[\ul{e}]$ and transverse to
$\ell$ then $D(\psi)$ has a quadratic singularity at $[\ul{e}]$; thus
$\wt{S}_q$ is Du Val at $t$ by Criterion~(\ref{duvalse}). Finally assume that
$V(B)$ is smooth at $[\ul{e}]$ and that $\ell$ is tangent to $V(B)$ at
$[\ul{e}]$. We notice that
\begin{equation}\label{alpiudue}
  (\ell\cdot V(B))_{[\ul{e}]} = 2.
\end{equation}
In fact if this does not hold then $\ell\subset
V(B)$ because $V(B)$ is a quadric and hence
$\ell\cap V(C)\subset V(F,A,B,C)$; this
contradicts~(\ref{intvuota}). Let $0\le i\le 3$
be such that $e_i\not=0$ and let
$a,b,c\in\CC[\PP^3\sm V(X_i)]$ be the regular
functions  $a:=A/X_i$, $b:=B/X^2_i$,
$c:=C/X^3_i$. From~(\ref{alpiudue}) we get that
there exists an open (in the classical topology)
 $U\subset Q$ containing $[\ul{e}]$ and analytic coordinates $(x,y)$ on $U$
 centered at $[\ul{e}]$ such that
\begin{equation}\label{birignao}
b|_U =y+x^2,\qquad I(\ell\cap U)=(y).
\end{equation}
Then $a|_U=\lambda y^2$ and $c|_U =\mu$ with
$\lambda,\mu\in\CC\{x,y\}$ units. Let
$\lambda\cdot\mu=\sum_{i,j}f_{i,j}x^i y^j$, where
$f_{i,j}\in\CC$. Then
\begin{equation}\label{equaloc}
  (b^2-4a\cdot c)|_U\equiv (1-4
  f_{0,0})y^2+2y(x^2-2f_{1,0}xy-2f_{0,1}y^2) \mod{(x,y)^4}.
\end{equation}
If $4f_{0,0}\not=1$ then $D(\psi)$ has a quadratic singularity at $[\ul{e}]$
and hence $\wt{S}_q$ is Du Val at $t$ by Criterion~(\ref{duvalse}). On the
other hand if $4 f_{0,0}=1$ then the
  multiplicity of $D(\psi)$ at $[\ul{e}]$ is $3$ and the strict transform of
  $D(\psi)$ under the blow-up of $Q$ at $[\ul{e}]$ intersects the exceptional
  divisor in at least $2$
  distinct points; thus Criterion~(\ref{duvalse}) applies again and we
  get that $\wt{S}_q$ is Du Val at $t$. This finishes the proof that if
  Item~(1) above holds then $\wt{S}_q$ is Du Val at $t$. Now we assume that
 Item~(2) holds, i.e.~that $Q$ is a cone with vertex $[\ul{e}]$ over a smooth conic.
Let $\rho\cl\wh{Q}\to Q$ be the blow-up of
$[\ul{e}]$ and $R$ be the exceptional divisor of
$\rho$.  Let $\wh{D}(\psi)\subset\wh{Q}$ be the
strict transform of $D(\psi)$. Since
$0=A([\ul{e}])=B([\ul{e}])$ and
$C([\ul{e}])\not=0$ we get that
\begin{equation}\label{decompongo}
  \rho^{*}D(\psi)=\wh{D}(\psi)+R.
\end{equation}
Thus $\rho^{*}D(\psi)$ is reduced, and there is a
unique square-root of
$\cO_{\wh{Q}}(\rho^{*}D(\psi))$, namely
$\rho^{*}\cO_Q(2)$; let $\vf\cl W\to \wh{Q}$ be
the corresponding normal double cover with branch
divisor $\rho^{*}D(\psi)$. We have a natural map
$\zeta\cl W\to \wt{S}_q$ which is an isomorphism
outside $t$ and such that
$\zeta^{-1}(t)=\vf^{-1}(R)$. Furthermore the
dualizing sheaf $\om_W$ is locally-free because
$W$ has hypersurface singularities and we have
\begin{equation}\label{ahiahi}
  \om_{W}\cong\zeta^{*}\om_{\wt{S}_q}.
\end{equation}
Thus it suffices to prove that $W$ has Du Val
singularities at all points of $\vf^{-1}(R)$.
Since $\rho^{*}D(\psi)$ is smooth at all points
of $R\sm (supp\wh{D}(\psi))$ we get that $W$ is
smooth at points of $\left(\vf^{-1}(R)\sm
\vf^{-1}(supp\wh{D}(\psi))\right)$. Let
$\wh{V}(A,F)\subset\wh{Q}$ be the strict
transform of $V(A,F)\subset Q$; we have
\begin{equation}\label{eccointer}
 R\cap (supp\wh{D}(\psi))=R\cap\wh{V}(A,F).
\end{equation}
Either $V(A,F)$ consists of two lines
$\ell_1,\ell_2$ or it is a \lq\lq double
line\rq\rq supported on the line $\ell$. In the
first case $R\cap\wh{V}(A,F)$ consists of two
points $r_1,r_2$. One easily checks that
$\rho^{*}D(\psi)$ has a quadratic singularity at
$r_1$ and at $r_2$; thus $W$ is Du Val at
$\vf^{-1}(r_1),\vf^{-1}(r_2)$ by
Criterion~(\ref{duvalse}). In the second case
$R\cap\wh{V}(A,F)$ consists of a single point
$r$: one easily checks that the multiplicity of
$\rho^{*}D(\psi)$ at $r$ is at most $3$ and that
if it is equal to $3$  then the strict transform
of $\rho^{*}D(\psi)$ under the blow-up of $r$
intersects the exceptional divisor in $2$
distinct points; thus $W$ is Du Val at
$\vf^{-1}(r)$ by Criterion~(\ref{duvalse}).
\end{proof}
Now we prove that $S_p$ has Du Val singularities.
Let $k$ be as in~(\ref{puntising}) and write
$sing(Y)=\{p,p_1,\ldots,p_k\}$. If $k=0$ then
$U_p=S_p$ and hence by the above proposition
$S_p$ has Du Val singularities. In order to prove
that $S_p$ has Du Val singularities when $k>0$ we
study the relation between $S_p$ and $S_{p_i}$.
Let $r_i:=span(p,p_i)$; thus $r_i\subset Y$. Let
$\Si(r_i)\subset \Gr(2,\PP^5)$ be the subset
parametrizing planes containing $r_i$. If
$[\Lambda]\in\Si(r_i)$ then $Y|_{\Lambda}$ is an
effective divisor because $Y$ does not contain
planes and we have
\begin{equation}\label{conica}
  Y|_{\Lambda}=r_i+c,\qquad c\in|\cO_{\Lambda}(2)|.
\end{equation}
Let $\G_i^0\subset \Si(r_i)$ be the subset
parametrizing planes $\Lambda$ such that the
conic $c$ of~(\ref{conica}) is reducible and
$r_i\not\subset supp(c)$. Let $\G_i\subset
\Gr(2,\PP^5)$ be the closure of $\G_i^0$. Let
$[\Lambda]\in\G_i^0$; since $Y$ is singular at
$p$ and at $p_i$ we must have $p,p_i\in supp(c)$
and hence there is a unique decomposition
$c=\ell+\ell'$ with $p\in \ell$ and
$p_i\in\ell'$. Thus we have regular maps
\begin{equation}\label{piconizero}
{\begin{matrix}
    \G_i^0 & \overset{\pi_i^0}{\lra} & S_p \\
   [\Lambda] & \mapsto & [\ell]
\end{matrix}}
\qquad\qquad
{\begin{matrix}
    \G_i^0 & \overset{\tau_i^0}{\lra} & S_{p_i} \\
   [\Lambda] & \mapsto & [\ell']
\end{matrix}}
\end{equation}
As is easily verified the above maps
extend to regular maps
\begin{equation}\label{piconitauconi}
  \pi_i\cl\G_i\to S_p,\qquad\qquad \tau_i\cl\G_i\to S_{p_i}.
\end{equation}
The fiber of $\pi_i$ over a point of
$S_p\sm\{[r_i]\}$ consists of a single point, and
the same holds for the fiber of $\tau_i$ over a
point of $S_{p_i}\sm\{[r_i]\}$. By Item~(2) of
Claim~(\ref{pocopoco}) we know that $S_p$ and
$S_{p_i}$ are normal and hence $\pi_i$ and
$\tau_i$ define isomorphisms
\begin{equation}\label{pitauiso}
  (\G_i\sm\pi_i^{-1}([r_i]))\overset{\sim}{\lra}
  S_p\sm\{[r_i]\},\qquad
  (\G_i\sm\tau_i^{-1}([r_i]))\overset{\sim}{\lra}
  S_{p_i}\sm\{[r_i]\}.
\end{equation}
In particular $\pi_i$ and $\tau_i$ are birational maps and hence $S_p$ is
birational to $S_{p_i}$.
\begin{prp}\label{bedda}
Keep assumptions and notation as above. The
embedding $\G_i\hra \Si(r_i)\cong\PP^3$ realizes
$\G_i$ as a quartic surface. Furthermore
$\pi_i^{-1}\left(U_p\cup\{[r_i]\}\right)$ is an
open subset of $\G_i$ with Du Val singularities.
\end{prp}
\begin{proof}
Over $\Si(r_i)$ we have a tautological family of
conics: the conic over $[\Lambda]$ is given by
the divisor $c$ appearing in~(\ref{conica}). Thus
we have a discriminant divisor
$\Delta_i\subset\Si(r_i)$ locally defined by the
determinant of a symmetric matrix defining the
family of conics. We have $\G_i\subset
supp(\Delta_i)$, however
$\G_i\not=supp(\Delta_i)$. In fact let
\begin{equation}\label{omegaidef}
  \Omega_i:=\{[\Lambda]\in\Si(r_i)|\quad
  Y|_{\Lambda}=2r_i+\ell,\quad \ell\in|\cO_{\Lambda}(1)|\}
\end{equation}
Clearly $\Omega_i\subset supp(\Delta_i)$ and
\begin{equation}\label{supporto}
  supp(\Delta_i)=\G_i\cup\Omega_i.
\end{equation}
A plane $\Lambda$  is parametrized by a point of
$\Si(r_i)$ if and only if it is tangent to $Y$ at
each point of $r_i$. Let
\begin{equation}\label{elleconidef}
  L_i:=\bigcap\limits_{y\in r_i} \Theta_y Y.
\end{equation}
Since $Y$ is singular at $p$ and $p_i$ but
$r_i\not\subset sing(Y)$ the linear space $L_i$
is a hyperplane. Thus $\Omega_i$ is a plane and
hence it is an irreducible component of
$supp(\Delta_i)$. Now we write out explicit
equations for $\Omega_i$, $\G_i$, etc. Let
\begin{equation}\label{ixezeta}
  \ul{X}:=X_0,\ldots,X_3,\qquad  \ul{Z}:=Z_0,Z_1.
\end{equation}
Choose projective coordinates $[\ul{X},\ul{Z}]$ on $\PP^5$
so that
\begin{equation}\label{puntifond}
p=[0,\ldots,0,1,0],\qquad p_i=[0,\ldots,0,1].
\end{equation}
Thus $r_i=V(\ul{X})$ and we have an obvious
identification $\Si(r_i)\cong\PP^3_{[\ul{X}]}$.
Since $r_i\subset Y$ we have $Y=V(\sum_j A_j
X_j)$ where $A_j\in\CC[\ul{X},\ul{Z}]$ is
homogeneous of degree $2$. Since $Y$ is singular
at $p$ and $p_i$ we have
$0=A_j(0,\ldots,0,1,0)=A_j(0,\ldots,0,1)$. Thus
\begin{equation}\label{eccoaconj}
 A_j=B_j+
 C_j Z_0+D_j Z_1+F_j Z_0 Z_1
\end{equation}
where $B_j,C_j,D_j,F_j\in\CC[\ul{X}]$ are
homogeneous of degrees $2$, $1$, $1$ and $0$
respectively. An easy computation gives that
\begin{equation}\label{omegaconi}
  \Omega_i=V(\sum_j F_j X_j).
\end{equation}
Let $[\ul{X}]$ correspond to the plane $\Lambda$;
a straightforward computation gives that the
conic $c$ appearing in~(\ref{conica}) is defined
by the $3\tm 3$ symmetric matrix
\begin{equation}\label{deltaequaz}
  M_i:=\begin{pmatrix}
    \sum_j B_j X_j & \sum_j C_j X_j
    & \sum_j D_j X_j \\
    \sum_j C_j X_j & 0 & \sum_j F_j X_j \\
    \sum_j D_j X_j  & \sum_j F_j X_j  & 0 \
  \end{pmatrix}
\end{equation}
computed at $\ul{X}$. In particular we get that
\begin{equation}\label{nocommon}
  V\left(\sum_j B_j X_j,\sum_j C_j X_j,
 \sum_j D_j X_j, \sum_j F_j X_j\right)=\es.
\end{equation}
The divisor $\Delta_i$ is defined by
\begin{equation}\label{determatrix}
  \det M_i=\left(\sum_j F_j X_j\right)\cdot
 \left(\sum_{j,h} \left(2 C_j X_j D_h X_h
 -B_j X_j F_h X_h\right)\right).
\end{equation}
Let $P_i\in\CC[\ul{X}]$ be the second factor
appearing in the right-hand side
of~(\ref{determatrix}). It follows
from~(\ref{omegaconi}) and~(\ref{nocommon}) that
$P_i$ does not vanish identically on $\Omega_i$;
thus by Equality~(\ref{supporto}) the zero-set of
$P_i$ is equal to $\G_i$. By Item~(2) of
Claim~(\ref{pocopoco}) we know that $\G_i$ is
irreducibile and hence we get that
\begin{equation}\label{marvin}
  (P_i)=m_i\G_i
\end{equation}
for some positive integer $m_i$. Let
\begin{equation}\label{sceltadia}
  [\ul{e}]\in
  V(\sum_j F_j X_j, \sum_j C_j X_j,
 \sum_j D_j X_j).
\end{equation}
Then
\begin{equation}\label{partialmod}
 P_i(\ul{e})=0,\qquad
 \frac{\partial P_i}{\partial X_s}(\ul{e})=
 -F_s\sum_j B_j(\ul{e}) e_j.
\end{equation}
Since $F_s\not=0$ for some $0\le s\le 3$ and since $\sum_j B_j(\ul{e})
e_j\not=0$ by~(\ref{nocommon}) we get that
\begin{equation}\label{bastaconti}
  \text{if~(\ref{sceltadia}) holds then $P_i(\ul{e})=0$ and
  $\nabla P_i(\ul{e})\not=0$ .}
\end{equation}
This proves that the $m_i$ appearing
in~(\ref{marvin}) is equal to $1$; since $\deg
P_i=4$ we get that $\G_i$ is a quartic, defined
by the vanishing of $P_i$. Let's show that
$\pi_i^{-1}\left(U_p\cup\{[r_i]\}\right)$ is an
open subset of $\G_i$ with Du Val singularities.
The subset $\left(U_p\cup\{[r_i]\}\right)\subset
S_p$ is open, see~(\ref{fuorising}), and hence
$\pi_i^{-1}\left(U_p\cup\{[r_i]\}\right)$ is
open. Next we notice that if $[\Lambda]\in\G_i$
and $\pi_i([\Lambda])=[r_u]$ with $u\not=i$ then
\begin{equation}\label{altramappa}
\tau_i([\Lambda])=[span(p_i,p_u)].
\end{equation}
In fact $Y|_{\Lambda}=r_i+r_u+\ell$ and since $Y$
is singular at $p_i$ and at $p_u$ we get that
$\ell=span(p_i,p_u)$. From~(\ref{altramappa}) we
get that
\begin{equation}\label{imagotau}
  \tau_i\left(\pi_i^{-1}\left(U_p\cup\{[r_i]\}\right)
  \right)=U_{p_i}\cup\{[r_i]\}.
\end{equation}
Let
$[\Lambda]\in\pi_i^{-1}\left(U_p\cup\{[r_i]\}\right)$.
By~(\ref{imagotau}) one of the following holds:
\begin{itemize}
  \item [(1)]
  $\pi_i([\Lambda])\in U_p$.
 \item [(2)]
$\tau_i([\Lambda])\in U_{p_i}$.
  \item [(3)]
 $[\Lambda]\in \pi_i^{-1}([r_i])\cap\tau_i^{-1}([r_i])$.
\end{itemize}
Suppose that~(1) holds. By~(\ref{pitauiso}) the
map $\pi_i$  is a local isomorphism onto $S_p$ in
a neighborhood of $[\Lambda]$. Applying
Proposition~(\ref{duralex}) with $q=p$ we get
that $\G_i$ is Du Val at $[\Lambda]$. If~(2)
holds a similar proof works: we apply
Proposition~(\ref{duralex}) with $q=p_i$. Finally
suppose that~(3) holds. We claim that
\begin{eqnarray}
\pi_i^{-1}([r_i]) & = V(\sum_j F_j X_j,
 \sum_j D_j X_j),\label{prima}\\
\tau_i^{-1}([r_i]) & = V(\sum_j F_j X_j,
 \sum_j C_j X_j)\label{seconda}.
\end{eqnarray}
In fact let $[\Lambda]\in \pi_i^{-1}([r_i])$ and
let $[\ul{X}]$ be its projective coordinates.
Since $\Lambda\cap Y=2r_i+\ell$ where
$p_i\in\ell$ we have $[\Lambda]\in\Omega_i$ and
$span(p_i,[\ul{X},0,0])\subset\PP(C_{p_i}Y)$.
This gives that $\pi_i^{-1}([r_i])$ consists of
those points of the right-hand side
of~(\ref{prima}) which are contained in $\G_i$.
Since $\G_i$ is the zero-locus of $P_i$ we get
that the right-hand side of~(\ref{prima}) is
contained in $\G_i$; this proves~(\ref{prima}).
Exchanging the r\^oles of $p$ and $p_i$ we get
Equation~(\ref{seconda}).
From~(\ref{prima})-(\ref{seconda}) we get that
\begin{equation}\label{piintertau}
  \pi_i^{-1}([r_i])\cap\tau_i^{-1}([r_i])=
  V(\sum_j F_j X_j, \sum_j C_j X_j,\sum_j D_j X_j).
\end{equation}
By~(\ref{bastaconti}) we get that $\G_i$ is smooth at every point of
$\pi_i^{-1}([r_i])\cap\tau_i^{-1}([r_i])$.
\end{proof}
Suppose that $k>0$ where $k$ is given
by~(\ref{puntising}): we prove that $S_p$ has Du
Val singularities. By Proposition~(\ref{duralex})
we know that $U_p$ has Du Val singularities. It
remains to show that $S_p$ has a Du Val
singularity at each $[r_i]$, where $1\le i\le k$.
Let $\ul{X},\ul{Z}$ be as in~(\ref{ixezeta}) and
assume that~(\ref{puntifond}) holds.  Projection
of $S_p$ from $[r_i]$ defines an embedding
$Bl_{[r_i]}(S_p)\hra\PP^3_{[\ul{X}]}$; the image
of this embedding is $\G_i$ and it gives an
identification of $\pi_i\cl\G_i\to S_p$ with the
blow-up of $[r_i]$. In particular since $\deg
S_p=6$ and $\deg\G_i=4$ we get that
  $mult_{[r_i]}S_p=2$. On the other hand $S_p$ has
embedding dimension $3$ at $[r_i]$ by
Claim~(\ref{pocopoco}) and hence we get that
\begin{equation}\label{tiraindietro}
  \om_{\G_i}=\pi_i^{*}(\om_{S_p}).
\end{equation}
Let $\rho_i\cl\wt{\G}_i\to\G_i$ be the minimal
desingularization of the singularities belonging
to $\pi_i^{-1}([r_i])$. By
Proposition~(\ref{bedda}) we know that $\G_i$ has
Du Val singularities along $\pi_i^{-1}([r_i])$
and hence
\begin{equation}\label{tiratira}
  \om_{\wt{\G}_i}=\rho_i^{*}(\om_{\G_i}).
\end{equation}
The regular map $\pi_i\circ\rho_i\cl\wt{\G}_i\to
S_p$ gives a desingularization of the singular
point $[r_i]$ and
by~(\ref{tiraindietro})-(\ref{tiratira}) we have
\begin{equation}\label{}
  \om_{\wt{\G}_i}=(\pi_i\circ\rho_i)^{*}(\om_{S_p}).
\end{equation}
This proves that $S_p$ has a Du Val singularity
at $[r_i]$.
\subsubsection{Proof of
Proposition~(\ref{eccemappa})}\label{firenze}
Let $S_p^{sm}\subset S_p$  be the smooth locus of
$S_p$. We have a cylinder map
\begin{equation}\label{cilindro}
  cyl\cl H_2(S_p^{sm};\ZZ)\to
  H^4(Bl_{S_p}\PP(\Theta_p\PP^5);\ZZ)
\end{equation}
defined as follows. Let
\begin{equation}\label{mappascoppio}
  \pi\cl Bl_{S_p}\PP(\Theta_p\PP^5)\to
\PP(\Theta_p\PP^5)
\end{equation}
be the blow-down map. Given a homology class
$\alpha\in H_2(S_p^{sm};\ZZ)$ represented by an
oriented closed smooth real surface $\Si\subset
S_p^{sm}$  the oriented smooth real $4$-fold
$\pi^{-1}\Si$ is in the smooth locus of
$Bl_{S_p}\PP(\Theta_p\PP^5)$, hence $\pi^{-1}\Si$
has a well-defined Poincar\'e dual class
$PD(\pi^{-1}\Si)\in
H^4(Bl_{S_p}\PP(\Theta_p\PP^5);\ZZ)$ independent
of the choice of representative $\Si$: we set
$cyl(\alpha):=PD(\pi^{-1}\Si)$.  Now let
$\ldots,R_i,\ldots$ be the irreducible components
of the desingularization map $\wt{S}_p\to S_p$;
thus we have
\begin{equation}\label{inclujei}
  j\cl S_p^{sm}\hra\wt{S}_p,
  \qquad j(S_p^{sm})=
  \left(\wt{S}_p\sm\bigcup_i R_i\right).
\end{equation}
Since $S_p$ has du Val singularities the map
$H_2(j)$ is injective and it gives an
identification
\begin{equation}\label{eccolambda}
  H_2(S_p^{sm})=\{\alpha\in H_2(\wt{S}_p;\ZZ)|
  \ \la\alpha,R_i \ra=0\ \forall
  R_i\},
\end{equation}
where $\la\cdot,\cdot\ra$ is the intersection
pairing on $H_2(\wt{S}_p;\ZZ)$. If $\alpha\in
H_2(\wt{S}_p;\ZZ)$ is Poincar\'e dual to a class
in $T(\wt{S}_p)$ then $\alpha$ belongs to the
right-hand side of~(\ref{eccolambda}). Thus via
Poincar\'e duality we get an injection
\begin{equation}\label{trascinu}
  T(\wt{S}_p)\hra H_2(S_p^{sm};\ZZ).
\end{equation}
Composing the above inclusion with the cylinder
map~(\ref{cilindro}) and tensoring with $\CC$ we
get a map
\begin{equation}\label{tsunami}
  \wt{\g}\cl T(\wt{S}_p)_\CC\lra
  H^4(Bl_{S_p}\PP(\Theta_p\PP^5)).
\end{equation}
A moment's thought will convince the reader that
the map above is a morphism of type $(1,1)$ of
Hodge structures. Furthermore for
$\alpha,\beta\in T(\wt{S}_p)_\CC$ we have
\begin{equation}\label{menoint}
  \int_{Bl_{S_p}\PP(\Theta_p\PP^5)}\wt{\g}(\alpha)\wedge
  \wt{\g}(\beta)=-\int_{\wt{S}_p}\alpha\wedge\beta.
\end{equation}
In fact this follows from a standard computation
 based on the fact that the normal bundle of the
 exceptional divisor of~(\ref{mappascoppio})
  has degree $-1$ on a
   fiber of the $\PP^1$-bundle $\pi^{-1}(S_p)\to
   S_p$. By Isomorphism~(\ref{duescoppi}) we may replace
the right-hand side of~(\ref{tsunami}) by
$H^4(Bl_p Y)$; thus $\wt{\g}$ defines a morphism
(we do not change its name) of type $(1,1)$
\begin{equation}\label{allafine}
  \wt{\g}\cl T(\wt{S}_p)_\CC\lra
  H^4(Bl_p Y).
\end{equation}
Let
\begin{equation}\label{blowblow}
  \rho\cl Bl_p Y\to Y
\end{equation}
be the blow-down map. The exceptional divisor of $\rho$ is the projectivized
normal cone $\PP(C_p Y)$. Composing the map of~(\ref{allafine}) with the
restriction map $H^4(Bl_p Y)\to H^4(\PP(C_p Y))$ we get
\begin{equation}\label{melampo}
  T(\wt{S}_p)_\CC\to H^4(\PP(C_p Y)).
\end{equation}
We claim that the above map is zero. It suffices to prove triviality of the
map
\begin{equation}\label{chiusuralampo}
  T(\wt{S}_p)_\CC\to H^4(\PP(C_p Y))/W_3 H^4(\PP(C_p Y))
\end{equation}
obtained by composing~(\ref{melampo}) with the
quotient map. The right-hand side
of~(\ref{chiusuralampo}) is a sub Hodge structure
of $H^4$ of any desing ularization of $\PP(C_p
Y)$; since $\PP(C_p Y)$ is a quadric we get that
the right-hand side of~(\ref{chiusuralampo}) is
of pure type $(2,2)$. By~(\ref{accaduez}) we get
that~(\ref{chiusuralampo}) has a non-zero kernel,
and since $T(\wt{S}_p)_\CC$ has no non-trivial
rational sub-Hodge structure we get that the
kernel of~(\ref{chiusuralampo}) is all of
$T(\wt{S}_p)_\CC$. Thus~(\ref{melampo}) is zero
and $Im(\wt{\g})\subset Im H^4(\rho)$ where
$\rho$ is the blow-down map~(\ref{blowblow}).
Hence there exists a morphism of type $(1,1)$ of
Hodge structures
\begin{equation}\label{ultimo}
  \widehat{\g}\cl T(\wt{S}_p)_\CC\lra
  H^4(Y)/\ker(\rho^{*}).
\end{equation}
such that  $\wt{\g}=H^4(\rho)\circ \widehat{\g}$.
Clearly $\ker(\rho^{*})\subset W_3 H^4(Y)$; we
let $\g$ be the composition of $\widehat{\g}$
with the quotient map $H^4(Y)/\ker(\rho^{*})\to
Gr^W_4 H^4(Y)$. This defines the morphism of
Hodge structures~(\ref{unouno}).
Equation~(\ref{chiave}) follows from
Equation~(\ref{menoint}).
\subsubsection{Comments}
Following is an example of $X$ a numerical $(K3)^{[2]}$ and $H$ a big and nef
divisor on $X$ with $(c_1(H),c_1(H))=2$ such that $f\cl X\to |H|^{\vee}$ is a
regular double covering of a cubic hypersurface - we do not know of any such
example with $H$ ample. Let $V$ be a $3$-dimensional complex vector space and
$\pi\cl S\to\PP(V)$ be a double covering ramified over a smooth sextic curve;
thus $S$ is a $K3$ surface. Let $X:=S^{[2]}$ and let $f$ be the composition
\begin{equation}\label{}
 S^{[2]}\to S^{(2)}\to \PP(V)^{(2)}\hra \PP(Sym^2 V)\cong\PP^5.
\end{equation}
The image of $\PP(V)^{(2)}\hra \PP(Sym^2 V)$ is the discriminant cubic
hypersurface; since $f$ has degree $4$ onto its image we get that $\int_X
c_1(H)^4=12$ and hence $(c_1(H),c_1(H))=2$ by~(\ref{hilbcase}). The divisor $H$
is big and nef and $f$ can be identified with the natural map $f\cl X\to
|H|^{\vee}$: thus $f$ has the stated properties.
\subsection{(5) of
Proposition~(\ref{varicasi}) does not hold}
 \setcounter{equation}{0}
In Subsubsection~(\ref{virus}) we will prove the
following result.
\begin{prp}\label{pieceofcake}
Let $Y\subset\PP^5$ be a quartic hypersurface
such that $\dim(sing Y)\ge 3$. Then $Y$ contains
a plane.
\end{prp}
Granting the above proposition  let's prove that
Item~(5) of Proposition~(\ref{varicasi}) does not
hold. We argue by contradiction. Assume that we
have $f\cl X\to Y$ regular of degree $3$ onto a
quartic hypersurface $Y\subset\PP^5$. By
Propositions~(\ref{nonredux})
and~(\ref{pieceofcake}) we get that $\dim(sing
Y)\le 2$. Let $R\in Div(X)$ be the ramification
divisor of $f$. Applying the adjunction formula
to $Y^{sm}:=(Y\sm sing Y)$ and Hurwitz' formula
to $f^{-1}(Y^{sm})\overset{f}{\to} Y^{sm}$ we get
that
\begin{equation}\label{ramificazione}
R\in|\cO_X(2H)|.
\end{equation}
By applying~(\ref{explicitrr}) we get that
\begin{equation}\label{propriouguali}
h^0(\cO_X(2H))=21=h^0(\cO_Y(2)).
\end{equation}
Thus the pull-back map $f^{*}\cl H^0(\cO_Y(2))\to
H^0(\cO_X(2H))$ is  an isomorphism and
from~(\ref{ramificazione}) we get that there
exists an effective Cartier divisor $D\in Div(Y)$
such that $f^{*}D=R$. Comparing the orders of
vanishing of $f^{*}D$ and $R$ at a prime
component of $R$ we get a contradiction.
\subsubsection{Proof of
Proposition~(\ref{pieceofcake})}\label{virus}
If $Y$ is not reduced or not irreducible then there is an irreducible component
of $Y$ of degree at most $2$ and the result follows immediately. Thus we may
assume that $Y$ is irreducible and reduced. Let $V$ be an irreducible component
of $sing Y$; intersecting $Y$ with a generic plane we get that $\deg V\le 3$. If
$\deg V=1$ there is nothing to prove. Assume that $\deg V=2$. If $V$ is singular
then $V$ contains planes and we are done. Thus we may assume that $V$ is smooth.
Let $L:=span(V)$. Then $L\cong\PP^4$ and $V$ is a quadric hypersurface in $L$.
Since $Y$ is irreducible of degree $4$ we have the cycle-theoretic intersection
\begin{equation}\label{interset}
  Y\cdot L=2V.
\end{equation}
We claim that there exists a complete
intersection of two quadrics
\begin{equation}\label{duequadriche}
  \wt{Y}=Q_1\cap Q_2\subset\PP^6
\end{equation}
such that $Y$ is isomorphic to the projection of
$\wt{Y}$ from a point outside $\wt{Y}$. In fact
let $\cI_V\subset\cO_{\PP^5}$ be the ideal sheaf
of (the reduced) $V$. The linear system
$|\cI_V(2)|$ has dimension $6$.  The rational map
\begin{equation}\label{ravel}
  \vf\cl\PP^5\cdots>|\cI_V(2)|^{\vee}\cong\PP^6
\end{equation}
is the composition of the blow-up of $V$ and
contraction of the strict transform of $L$ to a
point, call it $p$. The image of $\vf$ is a
smooth quadric $Q_1\subset\PP^6$. The inverse of
$\PP^5\cdots> Q_1$ is projection from $p$. The
image (strict transform) of $Y$ under $\vf$ is a
codimension-$1$ subset $\wt{Y}\subset Q_1$ which
does not intersect $p$ - use~(\ref{interset}) to
get this last statement. Thus $\deg\wt{Y}=\deg
 Y=4$ and hence there exists a quadric $Q_2\subset\PP^6$
such that~(\ref{duequadriche}) holds. By a
theorem of Debarre-Manivel~\cite{debman} we get
that $\wt{Y}$ contains a plane $\Lambda$. Since
projection from $p$ will map $\Lambda$  to a
plane in $Y$ we are done. Finally assume that
$\deg V=3$. The variety is non-degenerate: in
fact if $\dim(span(V))=4$ then $span(V)\subset Y$
contradiction. Since $V$ is non-degenerate of
degree $3$ we get that $V$ is smooth and linearly
normal; as is well-known~\cite{ion} it follows
that $V$ is the Segre $3$-fold
i.e.~$\PP^1\tm\PP^2$ embedded by
$\cO_{\PP^1}(1)\boxtimes\cO_{\PP^2}(1)$. Since
the Segre $3$-fold contains planes we are done.
\vskip 1cm
 \scriptsize{
Kieran G. O'Grady\\
Universit\`a di Roma ``La Sapienza",\\
Dipartimento di Matematica ``Guido Castelnuovo",\\
Piazzale Aldo Moro n.~5, 00185 Rome, Italy,\\
e-mail: {\tt ogrady@mat.uniroma1.it}. }

\end{document}